\documentclass[10pt]{amsart}

\usepackage{amsmath, amsfonts, amsthm}
\usepackage{algorithm,algorithmic}
\usepackage{cite}
\usepackage{booktabs}
\usepackage{geometry}
\usepackage{graphicx}
\usepackage{algorithm}
\usepackage{algorithmic}
\usepackage{epstopdf}
\usepackage{color}
\usepackage{subfigure}
\usepackage{pifont}
\usepackage{amsthm}
\usepackage{enumitem}
\usepackage{etoolbox}

\newtheorem{thm}{Theorem}[section]

\newtheorem{rmk}{Remark}[section]
\newtheorem{dfi}{Definition}[section]

\let\bbordermatrix\bordermatrix
\patchcmd{\bbordermatrix}{8.75}{4.75}{}{}
\patchcmd{\bbordermatrix}{\left(}{\left[}{}{}
\patchcmd{\bbordermatrix}{\right)}{\right]}{}{}

\begin{document}

\title[The PSL deflation method in reservoir simulation]{The parallel subdomain-levelset deflation method in reservoir simulation}

\author{J.H. van der Linden$^{\dag}$}
\author{T.B. J\"onsth\"ovel$^{\ddag}$}
\author{A.A. Lukyanov$^{\sharp}$}
\author{C. Vuik$^{\P}$}

\keywords{deflation, harmonic Ritz deflation, physics-based deflation, preconditioners, reservoir simulation, extreme eigenvalues, GMRES}
\subjclass[2000]{65F10, 65F08, 65Z05}

\begin{abstract}
Extreme and isolated eigenvalues are known to be harmful to the convergence of an iterative solver. These eigenvalues can be produced by strong heterogeneity in the underlying physics. We can improve the quality of the spectrum by `deflating' the harmful eigenvalues. In this work, deflation is applied to linear systems in reservoir simulation. In particular, large, sudden differences in the permeability produce extreme eigenvalues. The number and magnitude of these eigenvalues is linked to the number and magnitude of the permeability jumps. Two deflation methods are discussed. Firstly, we state that harmonic Ritz eigenvector deflation, which computes the deflation vectors from the information produced by the linear solver, is unfeasible in modern reservoir simulation due to high costs and lack of parallelism. Secondly, we test a physics-based subdomain-levelset deflation algorithm that constructs the deflation vectors a priori. Numerical experiments show that both methods can improve the performance of the linear solver. We highlight the fact that subdomain-levelset deflation is particularly suitable for a parallel implementation. For cases with well-defined permeability jumps of a factor $10^4$ or higher, parallel physics-based deflation has potential in commercial applications. In particular, the good scalability of parallel subdomain-levelset deflation combined with the robust parallel preconditioner for 
deflated system suggests the use of this method as an alternative for AMG.
\end{abstract}

\thanks{$^{\dag}$Delft University of Technology, Faculty of Electrical Engineering, Mathematics and Computer Science, Delft Institute of Applied Mathematics, 2628CN Delft, the Netherlands (joosthvanderlinden@gmail.com, c.vuik@tudelft.nl).}
\thanks{$^{\ddag}$Schlumberger Abingdon Technology Center, OX14 1UJ Abingdon, United Kingdom (tjonsthovel@slb.com).}
\thanks{$^{\sharp}$Schlumberger-Doll Research, Cambridge, MA 02139, USA (alukyanov@slb.com).}

\maketitle

\section{Introduction}
Recent challenges in the petroleum industry include managing larger data sets, providing higher field resolutions and computing more accurate multiphase flow predictions. These challenges entail complex geometries and high physical contrasts in the geological formations of petroleum reservoirs. At the same time, advancements in hardware, such as (cheap) parallel systems and GPU-acceleration, demand the development of new algorithms that utilize hardware innovations to exceed previous performance records. The continuous interplay between computational demand and supply has fueled the development of advanced reservoir simulation software.

At the core of any reservoir simulator is the solver mechanism. Modern reservoir simulators typically employ the Newton-Raphson method to solve the non-linear governing equations for a given timestep. The corresponding Jacobian matrix and linear system are solved by the Flexible Generalized Minimum Residual Method (FGMRES) \cite{Saad86} preconditioned by the Constrained Pressure Residual (CPR) preconditioner \cite{Wallis83,Wallis85,Cao05}. CPR decouples the linear system into two sets of equations, exploiting the specific properties of the pressure equation and transport equations. The former is solved with an Algebraic Multigrid (AMG) preconditioner \cite{Ruge87}, while the fully coupled system is solved using an ILU preconditioner. In this paper, the potential of an alternative for AMG only preconditoner is investigated by combing deflation method with different precondtioners (e.g., Jacobi). By removing unfavorable eigenvalues from the spectrum of the linear system, deflation can be used to improve convergence.

AMG is currently an industry standard for solving elliptic or parabolic partial differential equations. The method is  robust and scalable for a fixed problem size per processor. For a fixed total problem size, however, AMG is difficult to scale. In practice, creating reliable simulations with an increased number of grid cells is expensive. Therefore, while the number of available processors increases, reservoir engineers will often work with existing
solution strategies due to scalability issues. Despite the fact that AMG is optimal for serial computations of 
the pressure equation, the lack of strong scalability fuels the continued interest in alternatives for AMG.
As a result, two-level multiscale solvers (MS) were developed over the past decade in order to construct an accurate coarse-scale system 
honoring the fine-scale heterogeneous data (See, e.g., \cite{Hou97,Jenny-et-al-03,Jenny14,EffendievBookSiam,imsfv-jcp,yixuan-ams}).
The multiscale coarse-scale system is governed on the basis of locally computed basis functions, subject to reduced-dimensional boundary 
conditions and zero right-hand-side (RHS) terms. Multiscale solvers are naturally scalable and can be used as preconditioner.
Combing this method with deflation strategy (which has a lower algebraic complexity and is inexpensive to set up) may lead 
to a robust alternative of AMG. The preferred method of deflation in this paper (subdomain-levelset deflation) is also 
(strongly) scalable which will be combined with Jacobi and AMG preconditioners. The combination of multiscale solver and deflation 
method is subject of future research. 

Deflation was first proposed for symmetric linear systems and the conjugate gradient method (see e.g. \cite{Hestenes52}) by Nicolaides \cite{Nicolaides87} and Dost\' al \cite{Dostal88}. Both construct a deflation subspace consisting of deflation vectors to deflate unfavorable eigenvalues from the linear system. A range of deflation algorithms have been developed since, differing primarily in the method of application of the deflation operator and the approach to construct the deflation vectors. Deflation has been used with excellent results in a large number of applications, including electromagnetics \cite{Gersem01}, bubbly flow \cite{Tang07a,Tang07b,Tang07c,MacLachlan08}, structural mechanics and composite materials \cite{Jonsthovel09,Jonsthovel11,Jonsthovel13,Lingen12}, unsteady turbulent airfoil problems \cite{Carpenter10} and wave models in ship simulations \cite{vantWout10}. The work by Vuik and co-authors on layered problems in reservoir simulation \cite{Vuik98,Vuik99,Vuik00,Vuik01,Vuik02} is the foundation for this paper.

In \cite{Klie07}, the authors interweave algebraic multigrid cycles with deflation to solve several cases characterized by high permeability contrasts. The results are encouraging, showing improved convergence rates. We will argue, however, that Harmonic Ritz deflation, as used in \cite{Klie07}, is unfeasible in commercial applications due to the number of iterations required to compute the deflation vectors. Another popular approach is to combine deflation with a preconditioner based on the partial solution in a two-stage method. In \cite{Aksoylu07}, the partial solution is obtained in high permeability regions. In a related approach, the authors in \cite{Klie09} solve in aggregates of nodes with similar connectivity strength. Similar to our experiences, these methods work best for high physical contrasts.

Central to the investigation by Vuik and co-authors is the relation between the occurrence of extreme eigenvalues and large jumps in the PDE coefficients. In \cite{Vuik99}, the number of extreme eigenvalues is proven to be equal to the number of high-permeability layers (e.g. sand) between low-permeability layers (e.g. shale) for the diagonally scaled system matrix. Having observed this, the question arises how to utilize the predictable spectrum in layered problems. In \cite{Vuik99} and subsequent work, it is shown that the subspace spanned by the eigenvectors corresponding to the extreme eigenvalues can be approximated by a pre-determined space of algebraic deflation vectors. Convergence of the deflated CG method is shown to be independent of the size of the jumps in the coefficients.

We extend the work by Vuik and co-authors to non-symmetric linear systems arising from the fluid flow in porous media. Our serial physics-based deflation approach is based on the levelset deflation method \cite{Tang07b}. Regions of approximately constant permeability, separated by large jumps, are identified, and used to construct the deflation vectors. These vectors prove to be good approximations of the eigenvectors corresponding to the extreme eigenvalues caused by the jumps. Moreover, we will argue that the levelset deflation method allows for an efficient parallel implementation. Used in parallel, our deflation algorithm becomes very similar to the subdomain-levelset deflation method \cite{Tang07b}. In parallel subdomain-levelset deflation, the levelset deflation method is applied to each parallel subdomain. We extend the work on parallelizing the subdomain deflation method in \cite{Frank01,Tang05} to parallelize the subdomain-levelset deflation method. We use numerical experiments for cases with varying size and degree of complexity to compare the performance of Harmonic-Ritz eigenvector deflation and subdomain-levelset deflation. 

In the first part of this paper we give a brief introduction to deflation theory. Subsequently we present and motivate the choices of the deflation vectors and provide several numerical experiments on real simulation cases. In the last part of this paper we summarize and discuss future work. 



\section{Reservoir simulation}
\label{prob_def}
Physical properties are captured in the coefficients of the reservoir equations. We mainly focus on the permeability, or, roughly, the ease with which a fluid can flow through the porous medium. The grid, coefficients and wells give rise to a coupled system of partial differential equations. For a simple derivation of the incompressible two-phase (oil and water) flow, we refer to \cite{Lee09}. Central to our discussion is the pressure equation, which is defined
for incompressible two-phase flow as
\begin{align}
\label{eq:pressure}
	-\nabla \cdot \lambda \nabla p = q,
\end{align}
where $p$ is the pressure, $\lambda = \lambda_w + \lambda_o$ is the total mobility, $\lambda_w$ is the water mobility, $\lambda_o$ is the oil mobility 
and $q = q_w + q_o$ is the summed contribution from sources and/or sinks. The mobility coefficient is computed by summing
\begin{align}
\label{eq:mobility}
	& \lambda_w = k(x) \frac{k_{r_w} (S_w)}{\mu_w},\ \mathrm{and} \\
	& \lambda_o = k(x) \frac{k_{r_o} (S_o)}{\mu_o}.
\end{align}
Here, $\mu_w$, $\mu_o$ are the water and oil phase viscosity respectively, $k$ is the absolute permeability and $k_{r_w}$, $k_{r_o}$ are the relative permeability of water and oil phase respectively, and $S_w$, $S_o$ are the water and oil phase saturation respectively. The absolute permeability depends on the geological structure of the reservoir, usually determined through some form of geophysical imaging and lab experiments. As the absolute permeability is determined a priori, we often refer to $k$ as the initial permeability when considering incompressible rock. The relative permeability is a function of the saturation and will vary throughout the simulation. Properties $k$, $k_{r_w}$, $k_{r_o}$, $\mu_w$, $\mu_o$ (or in general $\lambda_w$ and $\lambda_o$) may exhibit large jumps, although the initial permeability will often be the dominant coefficient \cite{Lee09}. Hence, the initial permeability (simply referred to as the permeability in the remainder of this paper) will be used in the deflation methods.

The reservoir equations are discretized using the upstream finite-volume method. Including the well equations, the coupled linear system can be expressed as
\begin{align}
	Ax = \begin{bmatrix} A_{rr} & A_{rw} \\ A_{wr} & A_{ww} \end{bmatrix} \begin{bmatrix} x_r \\ x_w \end{bmatrix} = \begin{bmatrix} b_r \\ b_w \end{bmatrix} = b.
\label{eq:coupledsys}
\end{align}
The subscript $r$ refers to reservoir and $w$ stands for well. We assume $A \in \mathbb{R}^{n \times n}$, $b_r, x_r \in \mathbb{R}^{nr}$, $b_w, x_w \in \mathbb{R}^{nw}$. The sub-matrices $A_{rr}$ and $A_{ww}$ are square $nr \times nr$ and $nw \times nw$ matrices respectively. Each sub-matrix of the Jacobian represents a derivative, e.g. $A_{wr}$ is the derivative matrix of the well equations with respect to the reservoir variables. 

\begin{figure}[H]
\centering
\includegraphics[width=5cm]{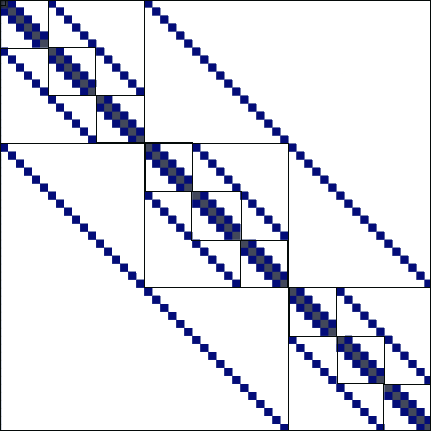}
\caption{$A_{rr}$ on a $6 \times 3 \times 3$ grid \cite{Schlumberger13}.}
\label{fig:ArrEx}
\end{figure}

\begin{figure}[H]
\centering
\includegraphics[width=7cm]{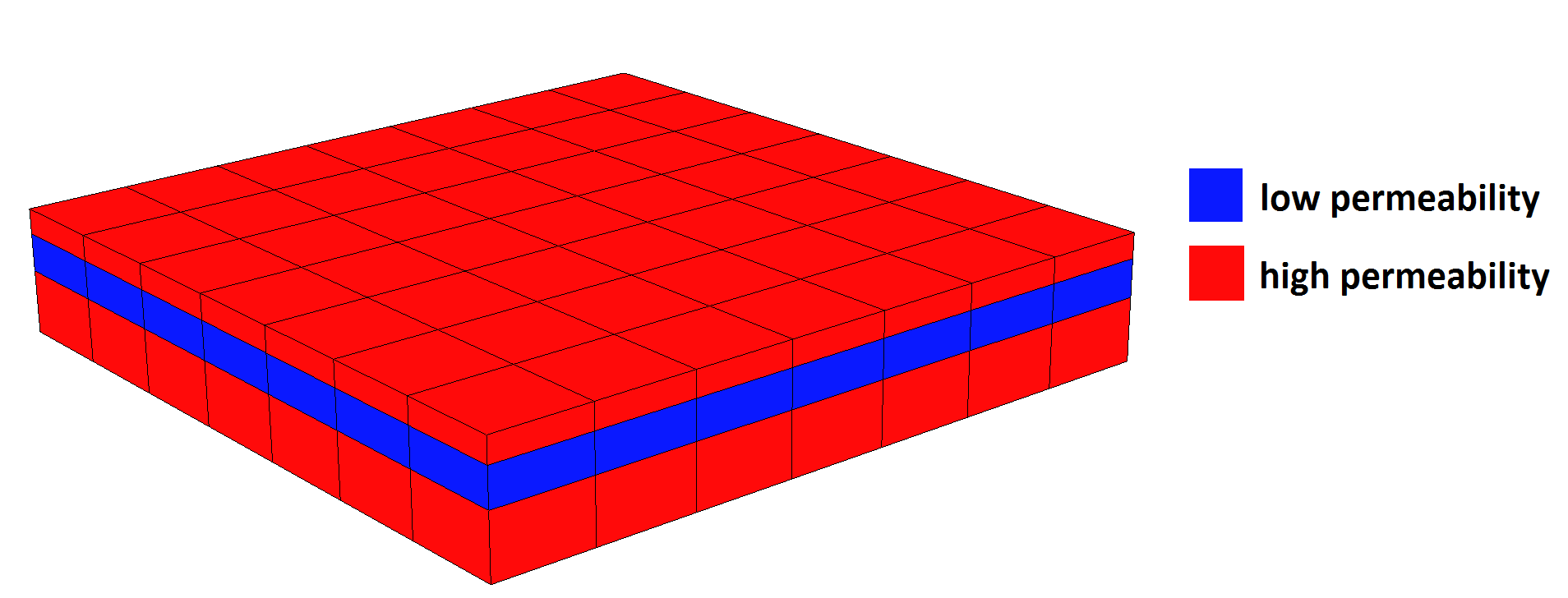}
\caption{Modified SPE5 permeability field.}
\label{fig:sandwich}
\end{figure}

\begin{figure}[H]
     \begin{center}
        \subfigure[] {
			\includegraphics[width=7cm]{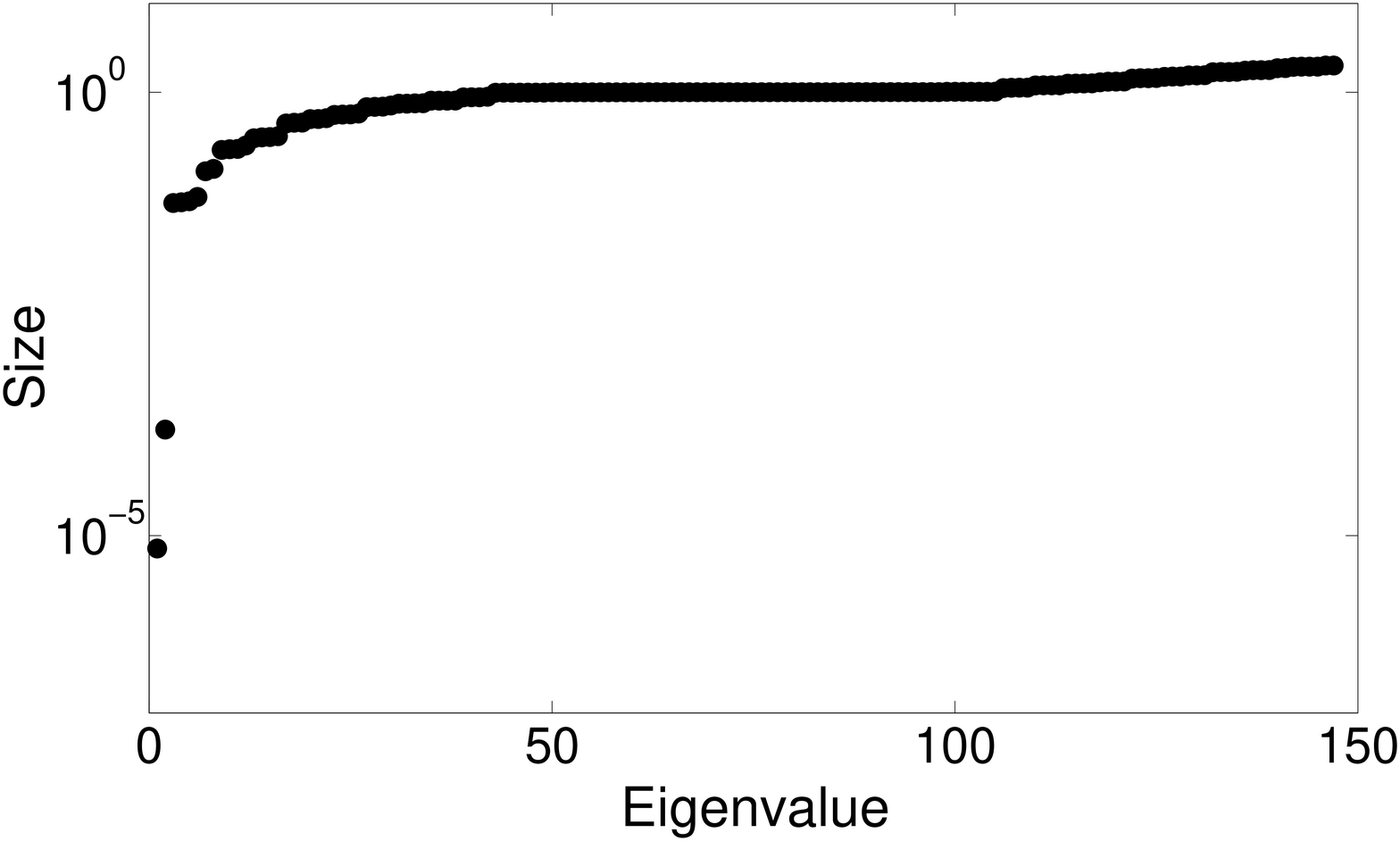}
				\label{fig:SPE5-sandwich-jump8-spectra} }
		\subfigure[] {
			\includegraphics[width=7cm]{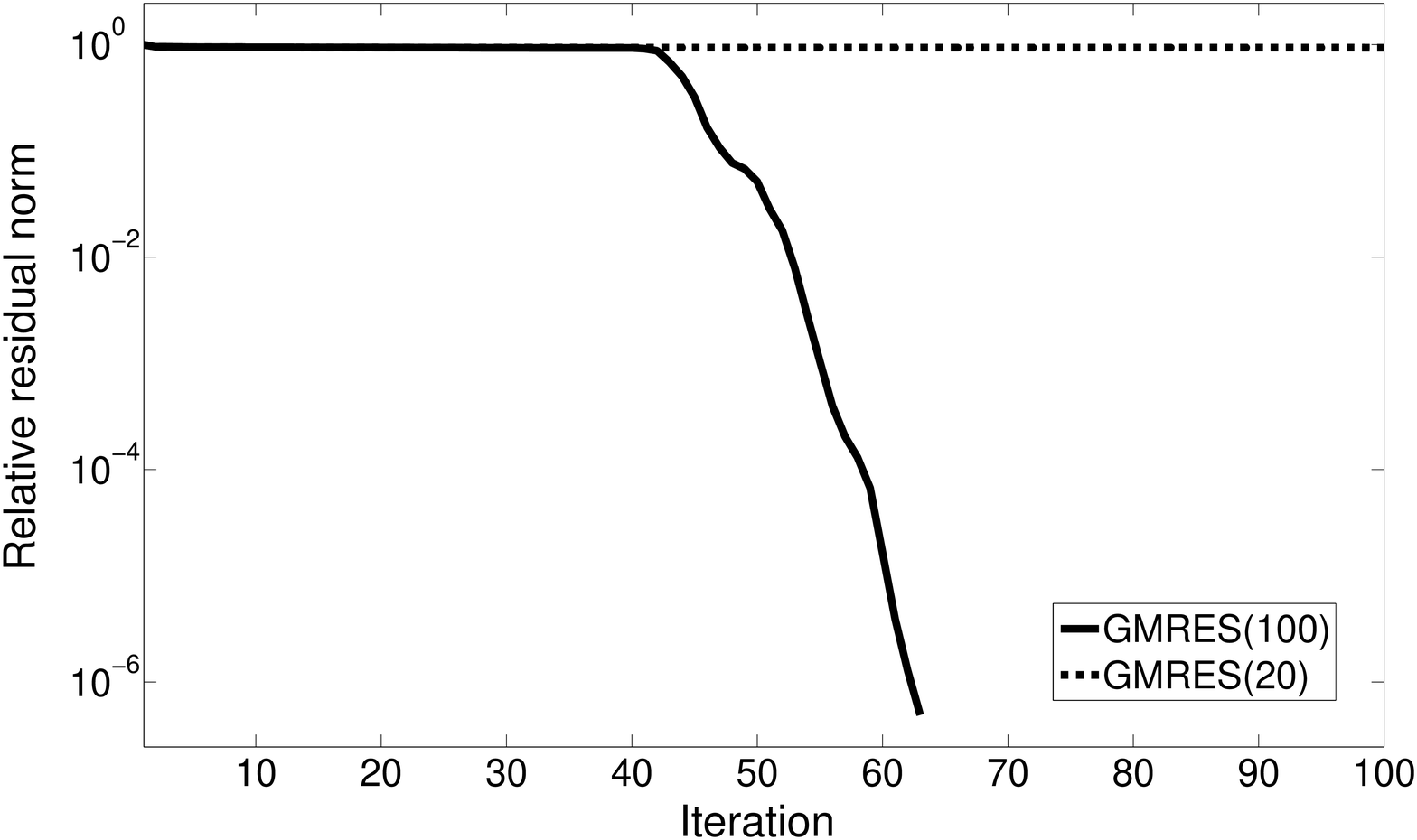}
				\label{fig:SPE5-sandwich-jump8-conv} }
	  \end{center}
	\caption{Pressure matrix spectrum of modified permeability field $\sigma = 10^6$ (a) and corresponding convergence of restarted GMRES (b).}
\label{fig:SPE5-sandwich-example}
\end{figure}

$A$ is typically very sparse. For example, in a three-dimensional $6 \times 3 \times 3$ grid, the matrix $A_{rr}$ is illustrated in Figure \ref{fig:ArrEx}. Squares indicate the block structure of $A$.
We assume that matrix $A$ has the following properties:
\begin{itemize}
	\item The computational grid consists of $n_c$ cells or `points', and each cell contains $n_u$ unknowns. Hence, $n = n_c \times n_u$.
	\item $A \neq A^T$, i.e. $A$ is non-symmetric.
	\item $\lambda \neq 0\ \forall\ \lambda \in \sigma (A)$, i.e. all eigenvalues are non-zero.
\end{itemize}
For illustrative purposes, consider the fifth comparative solution project of the Society of Petroleum Engineers \cite{Killough87} (denoted SPE5) with modified permeability, as shown in Figure \ref{fig:sandwich}.
A layer of low permeability is sandwiched between two layers of high permeability. The permeability is taken equal in $x$, $y$ and $z$ directions. The above pattern can be repeated to obtain a layered structure of alternating high/low permeability layers. We set the high permeability at $1$ and define $\epsilon$ as the order of the low permeability. For this permeability field, the following theorem is proven in  \cite{Vuik99}.
\begin{thm}[cf. {\cite[Theorem 3.1]{Vuik99}}]
For $\epsilon$ small enough, the diagonally scaled system matrix has only $n$ eigenvalues of $\mathcal{O}(\epsilon)$, where $n$ is the number of high-permeability layers lying between low-permeability layers.
\end{thm}
Here, homogeneous Neumann boundary conditions are used on all edges of the reservoir, except for the top boundary where Dirichlet boundary conditions are used. As expected, assigning value $1$ to the low-permeability layers and varying the value of the high-permeability layers (defined $\sigma$) yields a similar result. For example, for $\sigma = 10^6$, two extreme eigenvalues, corresponding to two jumps, are isolated from the main cluster in the spectrum of the pressure system. The result is illustrated in Figure \ref{fig:SPE5-sandwich-example}(a).

The intention in this paper is (a) to replace AMG with deflation as a preconditioner for the pressure system in (F)GMRES and (b) 
to combine deflation method with AMG where AMG is used as a preconditioner to the deflated system. The first task provides performance analysis of 
the deflation method whereas the second task sheds a light on overall performance of the deflated method with the robust preconditioner in a serial run.
As stated earlier, the combination of multiscale solver and deflation method is subject of future research. 

In Figure \ref{fig:SPE5-sandwich-example}(b), the convergence of solving the system giving rise to the spectrum in Figure \ref{fig:SPE5-sandwich-example}(a) using GMRES(20) and GMRES(100) is shown. Here, GMRES($m$) uses $m$ iterations before a restart.
Clearly, GMRES(20) has difficulty to converge, whereas GMRES(100) converges after about 40 iterations. In the next section, the concept of natural deflation is introduced to explain this behavior.

\section{Deflation}
The superlinear convergence of GMRES was associated by Van der Vorst and Vuik \cite{vanderVorst93} with the convergence of the Ritz values of the Hessenberg matrix to the eigenvalues of the operator $A$. If the Krylov subspace reaches a sufficient size, the Ritz values will be close to the eigenvalues of $A$. From that point on, GMRES will behave as if these approximated eigenvalues have been \textit{naturally deflated} from $A$, resulting in faster convergence. To illustrate this idea, the Ritz values corresponding to the convergence history in Figure \ref{fig:SPE5-sandwich-example}(b) are plotted below.
\begin{figure}[H]
     \begin{center}
        \subfigure[] {
            \includegraphics[width=7cm]{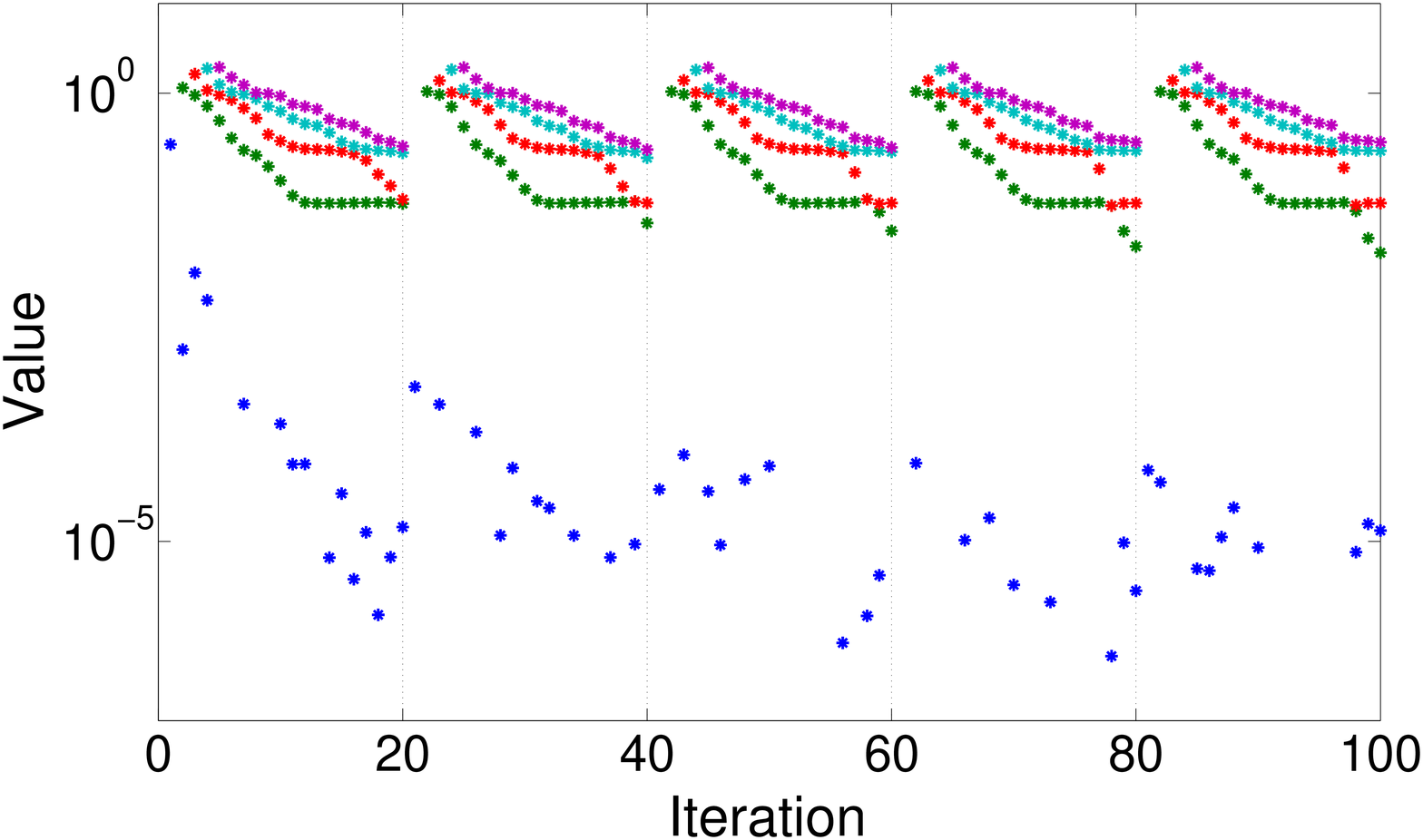}
						\label{fig:example_ritz_restarted} }
		\subfigure[] {
            \includegraphics[width=7cm]{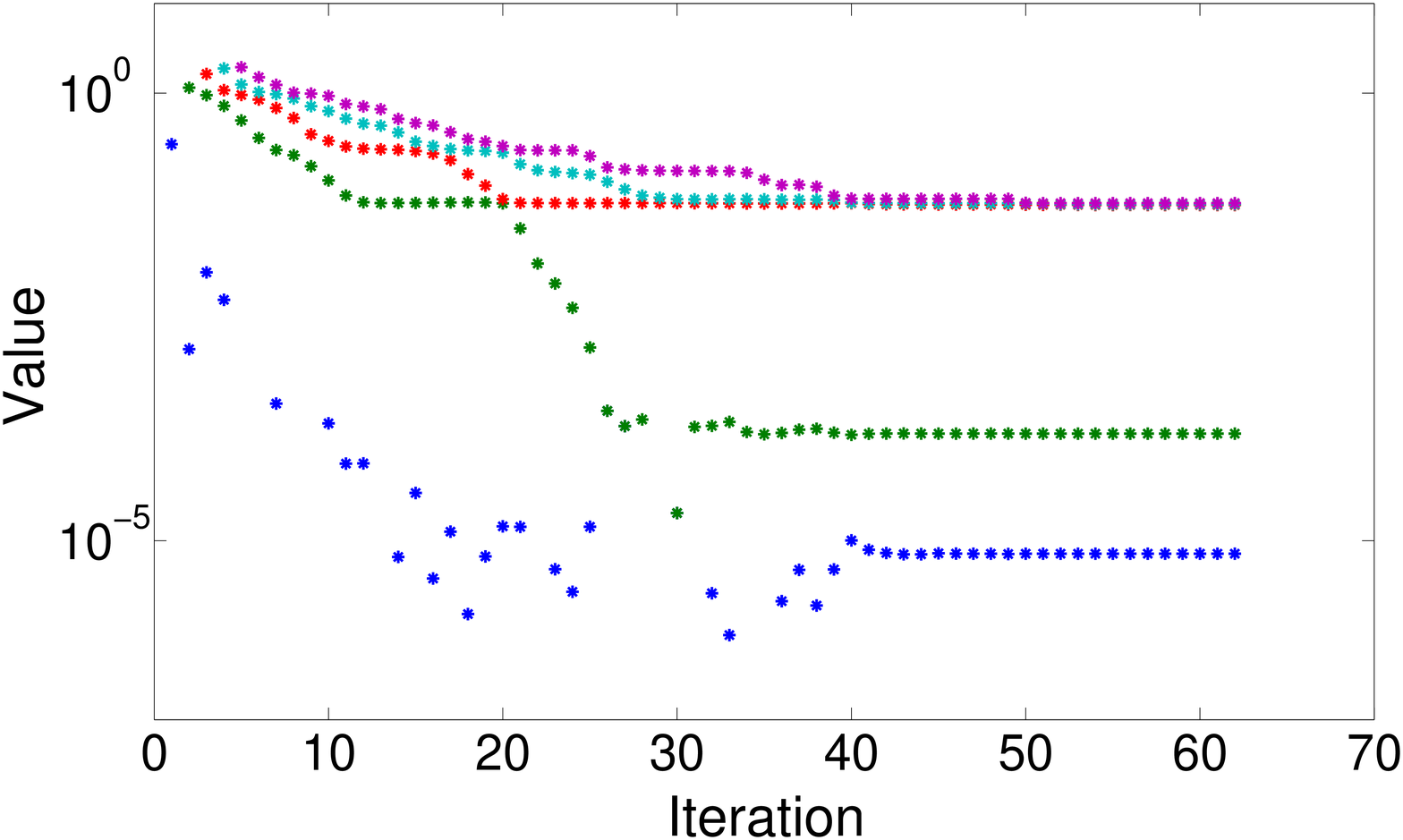}
						\label{fig:example_ritz_full} }					
    \end{center}
    \caption{Smallest Ritz values for (a) GMRES(20) and (b) GMRES(100).}
    \label{fig:example_ritz}
\end{figure}
In Figure \ref{fig:example_ritz}(a), the five smallest Ritz values of the GMRES(20) convergence in Figure \ref{fig:SPE5-sandwich-example}(b) are plotted. Clearly, the convergence of the Ritz values towards the exact eigenvalues is reset after each restart.  In Figure \ref{fig:example_ritz}(b), which corresponds to the GMRES(100) convergence, the two smallest Ritz values converge to the two smallest eigenvalues in Figure \ref{fig:SPE5-sandwich-example}(a). In conclusion, precisely as the Ritz values approach the exact eigenvalues, at about 40 iterations, convergence of GMRES(100) becomes superlinear.

Natural deflation occurs, because eigenvector components corresponding to eigenvalues are removed from the linear system. Extreme eigenvalues, in particular, are detrimental to the convergence of GMRES \cite{Frank01,Erhel96,Burrage98}. Unfortunately, a restart in GMRES erases the obtained Krylov subspace before it might reach sufficient size to allow for natural deflation, as illustrated in Figure \ref{fig:example_ritz}. In addition, even if the Krylov subspace grows big enough and deflation occurs, it will only aid to convergence in the current cycle. This discussion suggests convergence would improve if (small) eigenvalues could be removed \textit{artificially}. This idea gives rise to the method of deflation.

In the remainder of this section, we will give an overview of the fundamentals of the deflation theory, referring to relevant literature for more details. First, we will discuss the formulation and some results related to the convergence of GMRES. In the second half of the section, we will provide an overview of two methods to construct the deflation vectors.
\subsection{Deflation Theory}
For a thorough introduction of deflation methods for symmetric problems, we refer to \cite{Tang07}. We will follow the derivation for the non-symmetric case, discussed in \cite{Frank01}.
\begin{dfi}
\label{def:deflation}
Let $A \in \mathbb{R}^{n \times n}$ be the non-symmetric linear system matrix, and assume that the deflation matrix $Z \in \mathbb{R}^{n \times d}$ with $d$ deflation vectors is given. Then we define the matrix $E \in \mathbb{R}^{d \times d}$ as
\begin{align}
\label{eq:E}
	E = Z^T A Z,
\end{align} 
assuming $E^{-1}$ exists, and the deflation operators $P_1$ and $P_2$ as
\begin{align}
\label{eq:deflops}
\begin{split}
	P_1 &= I - A ZE^{-1}Z^T, \\
	P_2 &= I - ZE^{-1}Z^T A.
\end{split}
\end{align}
\end{dfi}
Inverting $E$ is relatively cheap, since, in general, $d << n$. In the symmetric case we have $P_1^T =  P_2$. If the columns of $Z$ (the deflation vectors) form an invariant subspace of $A$, then $P_1A$ has $d$ zero eigenvalues \cite{Yeung10}. In particular, when $Z$ contains $d$ exact eigenvectors of $A$, then applying $P_1$ to $A$ implies `deflating' the $d$ corresponding eigenvalues to zero.

$P_1$ and $P_2$ are used to apply deflation as follows. We split the solution $x$ of \eqref{eq:coupledsys} in two parts \cite{Frank01}:
\begin{align}
\label{eq:defsplit}
\begin{split}
	x &= (I - P_2)x + P_2 x \\
	  &= ZE^{-1}Z^T b + P_2 x,
\end{split}
\end{align}
which can be rewritten to obtain
\begin{align}
	P_1 A\hat{x} = P_1 b,
\label{eq:defsys}
\end{align}
where $x$ was replaced by $\hat{x}$ to denote the solution to the 'deflated system'. Note that $P_1 A$ has at least one zero eigenvalue, so the system is singular. Consequently, the solution $\hat{x}$ is not necessarily the solution of the original linear system $Ax = b$, as $\hat{x}$ may contain components in the nullspace of $P_1 A$. It can be easily shown that $P_2 \hat{x}= P_2 x$, however, see e.g. \cite[Lemma 3.5]{Tang07}. Hence, even though \eqref{eq:defsys} is singular, the projected solution $P_2 \hat{x}$ is unique because it has no components in the null space $\mathcal{N}(P_1 A)$. Note that 'deflated system' \eqref{eq:defsys} must be still efficiently solved to have overall robustness.

In conclusion, the solution to the original linear system can be found by solving \eqref{eq:defsys} for $\hat{x}$, and substituting the result in
\begin{align}
	x = ZE^{-1}Z^T b + P_2 \hat{x}.
\label{eq:projectback}
\end{align}
In the upcoming sections, we will use the notation $x^* = ZE^{-1}Z^Tb$. Note that if 
the columns of $Z$ (the deflation vectors) are multiscale basis functions, then 
$x^*$ is the solution similar to the one provided by algebraic multiscale method 
(AMS) \cite{yixuan-ams}.

\subsubsection{Convergence and condition number}
\label{sec:gmresconv}
Assume the columns of the deflation matrix $Z$ form a basis of some $A$-invariant subspace. Let $r_m$ and $\hat{r}_m$ be the $m$'th residual of the original linear system \eqref{eq:coupledsys} and the deflated system \eqref{eq:defsys}, respectively, solved using GMRES. Then, starting with the same initial guess, we have \cite[Theorem 5.1]{Yeung10}
\begin{align}
\label{eq:deflres}
	\|\hat{r}_m\|_2 \leq \|r_m\|_2 \qquad \forall\ m = 1,2,\ldots
\end{align}
In the symmetric case, the condition number of the deflated system is better than the condition number of the original linear system \cite[Theorem 2.2]{Frank01}. In the nonsymmetric case, this result does not hold, although similar convergence behavior was observed by the authors of \cite{Frank01} as long as the asymmetric part of $A$ is not too dominant.

In conclusion, as long as $Z$ contains a basis for an $A$-invariant subspace, deflated GMRES will converge faster than regular GMRES. Although a theoretical proof does not exist, we expect similar behavior when near-invariant subspaces are used. Since the eigenvectors are by construction an invariant subspace of $A$, the next section introduces methods to approximate the eigenspace.

\subsection{Deflation vectors}
A number of approaches to compute the deflation subspace $Z$ have been proposed in deflation-related literature, with varying degrees of effectiveness depending on the application. We will review the two most prominent methods. 

Since our linear system is non-symmetric, approximate eigenvalues and eigenvectors are either real or come in complex-conjugate pairs \cite[Theorem 1.3]{Stewart01}. To retain real arithmetic, a complex-conjugate pair of eigenvectors should be replaced by one eigenvector containing the real part of the complex pair, and one eigenvector containing the imaginary part. For a complex conjugate pair $(u_k,u_{k+1})$, this is done by using the transformation
\begin{align*}
	\begin{bmatrix} u'_k \\ u'_{k+1} \end{bmatrix} = \frac{1}{2} \begin{bmatrix} 1 & 1 \\ i & -i \end{bmatrix} \begin{bmatrix} u_k \\ u_{k+1}. \end{bmatrix}
\end{align*}
The vectors $(u'_k,u'_{k+1})$ are Schur vectors of $A$.

\subsubsection{Harmonic Ritz deflation}
\label{sec:harritzdefl}
For an approximate eigenvector $z$ of $A$ with corresponding eigenvalue $\theta$, the Galerkin orthogonal projection problem \cite{Chapman96} states
\begin{align*}
	Az - \theta z \perp \mathcal{K}_m,
\end{align*}
where $\mathcal{K}_m$ is the Krylov subspace. Using a basis $V_m$ of $\mathcal{K}_m$ with $z = V_m y$, this becomes
\begin{align}
	V_ m^T(A - \theta I)V_m y = 0.
\label{eq:oblique}
\end{align}
Using the well known identities $H_m = V_m^TAV_m$ and $V_m^T V_m = I$, equation \eqref{eq:oblique} reduces to
\begin{align*}
	H_m y = \theta y,\ z = V_m y.
\end{align*}
Ritz vectors approximate the eigenvectors of $A$. Moreover, the Ritz values tend to approximate the eigenvalues of $A$. Therefore, we can take the $d$ approximated eigenvectors $z$ corresponding to the $d$ smallest Ritz values as the columns of $Z$. In terms of eigenvector approximations for extreme eigenvalues, Chapman, Saad \cite{Chapman96} and Morgan \cite{Morgan02} report that Ritz vectors are outperformed by harmonic Ritz vectors. The latter concept will be introduced next.

Whereas Ritz vectors are formed by imposing a Galerkin projection, harmonic Ritz vectors are obtained by using the Petrov-Galerkin orthogonality conditions. The approximation error of the approximate eigenpair $(\theta, z)$ is set orthogonal to the subspace $A\mathcal{K}_m$, i.e.
\begin{align}
	Az - \theta z \perp A\mathcal{K}_m \qquad \Leftrightarrow \qquad (AV_m)^T(Az - \theta z) = 0,\ z = V_m y.
\label{eq:harritzcond}
\end{align}
There are generally two methods to solve \eqref{eq:harritzcond}. Denote $\lfloor B \rfloor$ as the matrix $B$ without its last row. Equation \eqref{eq:harritzcond} is equivalent to solving either:
\begin{align}
\label{eq:defleig}
\begin{split}
	&\mbox{(a) the eigenvalue problem } (H_m + h_{m+1,m}^2 H_m^{-T} e_m e_m^T ) y = \theta y\mbox{, or} \\
	&\mbox{(b) the generalized eigenvalue problem } R_m y = \theta \lfloor Q_m V_{m+1}^T V_m \rfloor y.
\end{split}
\end{align}
where $e_m$ is the $m$-th canonical vector in $\mathbb{R}^{m}$, the orthogonal matrix $Q_m$ and 
upper triangular $R_m$ matrix satisfy the condition $H^{T}_{m}=R_m\lfloor Q_m V_{m+1}^T V_m \rfloor$. Although both formulas result in a valid spectrum, approach (b) is preferable for computational reasons. Approach (a) would require keeping an original copy of the matrix $H_m$ in memory, without the Givens rotations. Approach (b), on the other hand, does not have this requirement, and allows the same Givens rotations $Q_m$ that are saved and applied to $H_m$, to be used on $V_{m+1}^T V_m$.

Note that as the Ritz values converge to the eigenvalues (e.g. Figure \ref{fig:example_ritz}), the harmonic Ritz eigenvector approximations approach the true eigenvectors. Therefore, the cycle size $m$ needs to be chosen sufficiently large in order to obtain reasonable approximations. In the results, we vary the value of $m$ in deflated GMRES using harmonic Ritz vectors, and demonstrate the impact on the convergence. Although $Z$ is not sparse in this case, the harmonic Ritz vectors can be computed at a relatively small cost. Since the user only has to specify how many vectors should be included in the deflation operator, the method has a black-box nature. 

The pseudocode for GMRES using harmonic Ritz deflation is given in Algorithm \ref{alg:gmres-defl-HR}.
\begin{algorithm}[H]
\caption{right-preconditioned GMRES with harmonic Ritz deflation}
\begin{algorithmic}[1]
\label{alg:gmres-defl-HR}
	\STATE Setup $P_1 = P_2 = I$, $x^* = 0$ and flag = false
	\STATE Compute $r_0 = P_1(b - Ax_0)$, $\beta = \|r_0\|_2$, and $v_1 = r_0 / \beta$.\label{alg:gmres-defl-dyn:line:start}
	\vspace{0.15cm}
	\FOR{$j=1,2,\ldots,m$}
		\STATE $w_j = P_1 AM^{-1}v_j$
		\FOR{$i=1,\ldots,j$}
			\STATE $h_{i,j} = (w_j,v_i)$
			\STATE $w_j = w_j - h_{ij}v_i$
		\ENDFOR
		\STATE $h_{j+1,j} = \|w_j\|_2$
		\IF{$h_{j+1,j} = 0$ \OR converged} 
			\STATE set $m = j$ and \textbf{go to} \ref{alg:gmres-defl-dyn:line:Hm}
		\ENDIF
		\STATE $v_{j+1} = w_j / h_{j+1,j}$
	\ENDFOR
	\vspace{0.15cm}
	\IF{flag = false} \label{alg:gmres-defl-dyn:line:ifstart} 
		\STATE solve $(AV_m)^T(Ay_k - \theta_k y_k) = 0$ for $y_k$, for $k = 1,\ldots,d$ \label{alg:gmres-defl-dyn:line:Ritzsolve} 
		\STATE $z_k = V_m y_k$, for $k = 1,\ldots,d$
		\STATE fill $Z = [z_1 \ldots z_d]$ 
		\STATE $A_p = AM^{-1}$, $E = Z^T A_p Z$
		\STATE $P_1 = I - A_p ZE^{-1}Z^T$, $P_2 = I - ZE^{-1}Z^T A_p$
		\STATE $x^* = Z E^{-1} Z^T b$
		\STATE set flag = true
	\ENDIF \label{alg:gmres-defl-dyn:line:ifend} 
	\vspace{0.15cm}
	\STATE Fill $\bar{H}_m = \{h_{ij}\}$ for $1 \leq i \leq m+1$, $1 \leq j \leq m$. \label{alg:gmres-defl-dyn:line:Hm}
	\STATE Compute the minimizer $u_m$ of $\|\beta e_1 - \bar{H}_m u\|_2$ and set $x_m = x_0 + M^{-1}V_m u_m$. 
	\STATE \algorithmicif\ converged \algorithmicthen\ $x_m = P_2x_m + x^*$ and return \algorithmicelse\ set $x_0 = x_m$ and \textbf{go to} 		
																												\ref{alg:gmres-defl-dyn:line:start} 
\end{algorithmic}
\end{algorithm}
Algorithm \ref{alg:gmres-defl-HR} assumes that the deflation vectors are not available at the start of the iteration. Instead, the information generated by GMRES is used to compute the harmonic Ritz vectors. $Z$ is constructed by taking the first $d$ eigenvector approximations $y_k$ in line \ref{alg:gmres-defl-dyn:line:Ritzsolve}. The generalized eigenvalue problem is solved as \ref{eq:defleig}(b). In the implementation we include a flag that freezes the deflation operator after the first restart. In theory, the harmonic Ritz vectors could be recomputed after every cycle, and appended to $Z$. This will further improve convergence, as more eigenvalues are deflated from the spectrum, but has two disadvantages. Firstly, the memory requirements increase as the size of $Z$ increases. Secondly, the computational costs of the matrix-vector products, inner products and Galerkin solve (in applying $P_1$ and $P_2$) increase as well. These costs, in practice, do not outweigh the time gain from the reduced number of iterations.\\
\begin{rmk}
\label{rmk:AM}
To apply $P_1$ and $P_2$ to a vector (these matrices are never explicitly formed) the use of $M^{-1}$ in $A_p = AM^{-1}$ requires two additional matrix-vector products with the preconditioner. The costs associated with these products would immediately render deflation unfeasible in practice. Therefore, we experimented with $A_p = A$ instead, and found that for right-preconditioning deflation will still perform well. This does not seem to be the case for left-preconditioning. The mathematical details are beyond the scope of this paper, but we plan on further analyzing this discrepancy in future work.
\end{rmk}

\subsubsection{Physics-based deflation}
\label{sec:domdefl}
Subdomain deflation has been introduced by Nicolaides \cite{Nicolaides87} and Mansfield \cite{Mansfield90,Mansfield91}. Let $\Omega$ be the computational domain, which is divided into $d$ nonoverlapping subdomains $\Omega_j$, $j = 1,\ldots,d$. After discretization, denoted by subscript $h$, let $x_i$ be a grid point in the discretized domain $\Omega_{h_j}$. We define the deflation vector $z_j$ corresponding to $\Omega_{h_j}$ as
\begin{align}
\label{eq:domdefl}
	(z_j)_i = \begin{cases} 1, \qquad x_i \in \Omega_{h_j} \\ 0, \qquad x_i \in \Omega_h \setminus \bar{\Omega}_{h_j}. \end{cases}
\end{align}
The deflation subspace is defined as $Z = [z_1, \ldots, z_d]$. The vectors in $Z$ are piecewise-constant, disjoint and orthogonal. For this choice of the deflation subspace, the deflation projectors $P_1$ and $P_2$ essentially aggregate each subdomain in a single cell. Hence, subdomain deflation is closely related to domain-decomposition methods and multigrid \cite{Frank01}. For problems in bubbly flow, the span of the deflation vectors \eqref{eq:domdefl} approximates the span of the eigenvectors corresponding to the smallest eigenvalues \cite{Tang07b}. Furthermore, subdomain deflation can be associated with fluid in place data (i.e., mobility distribution), which would be layered or channelized, in the immiscible black-oil case, commonly used in reservoir simulations.

In \cite{Vuik02}, a time-dependent diffusion equation is investigated for a layered medium representing the earth's crust. Three approaches are used to construct the domain-based deflation vectors. First, the authors require $(z_j)$ to satisfy the finite element discretization of the governing equation on all subdomains with low permeability. The deflation vectors satisfy \eqref{eq:domdefl} for the remaining highly permeable subdomains. This approach is robust for all test problems, yet costly due to the extra solves required. Second, the authors use the vectors \eqref{eq:domdefl} only on the high-permeability layers, and last, \eqref{eq:domdefl} is used for both the high- and low-permeability layers. The latter method turns out to be the most efficient and robust, and will be used in this work.

Vermolen et al. use subdomain deflation in \cite{Vermolen04} to solve a Poisson problem with discontinuous coefficients. In particular, the amount of overlap between subdomains is investigated. For large contrasts in the coefficients, the authors conclude that no overlap is the best choice. For no contrasts, on the other hand, average overlap is superior. This observation gives rise to the so called `weighted overlap' method, which mimics average and no overlap in the case of no contrasts and large contrasts, respectively. We obtained good results without overlap, but plan on experimenting with the weighted overlap approach to account for the mixed high/low permeability contrast that is often found in real reservoirs.

\begin{figure}[H]
     \begin{center}
        \subfigure[] {
            \includegraphics[width=4.5cm]{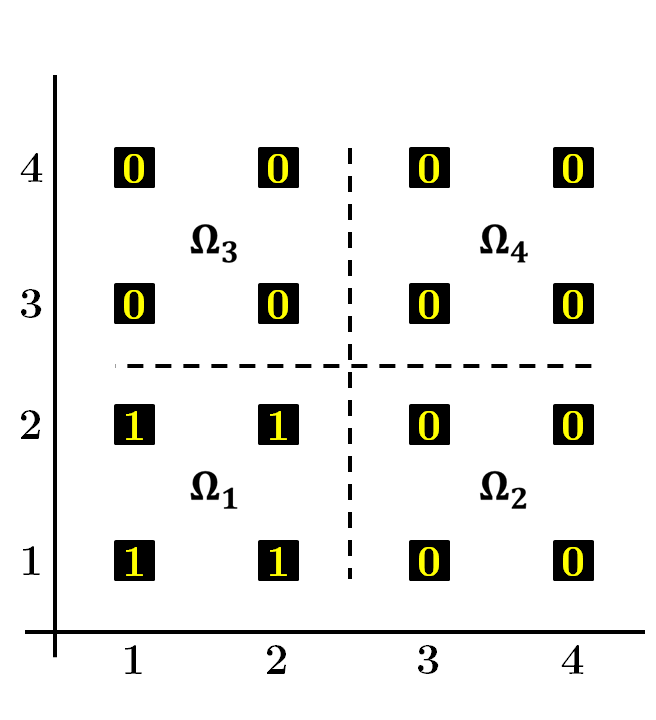}
						\label{fig:subdomain_deflation} }
		\subfigure[] {
            \includegraphics[width=4.5cm]{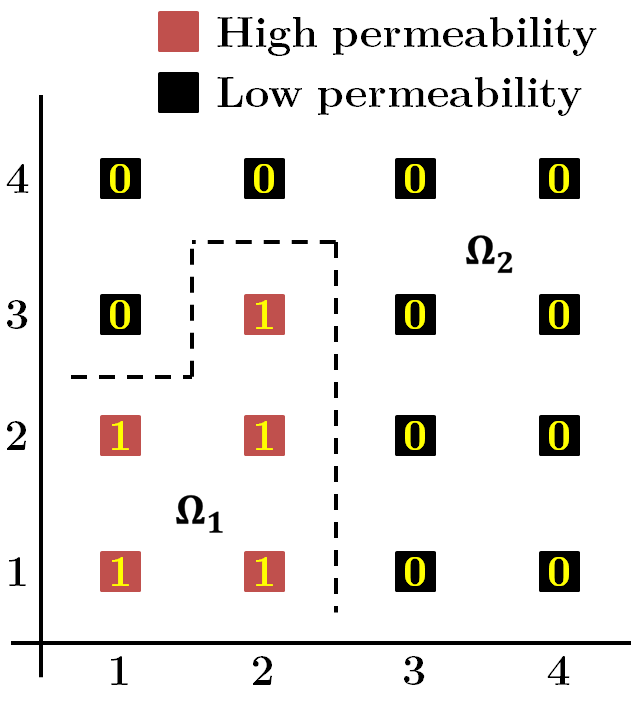}
						\label{fig:levelset_deflation} }	
		\subfigure[] {
            \includegraphics[width=4.5cm]{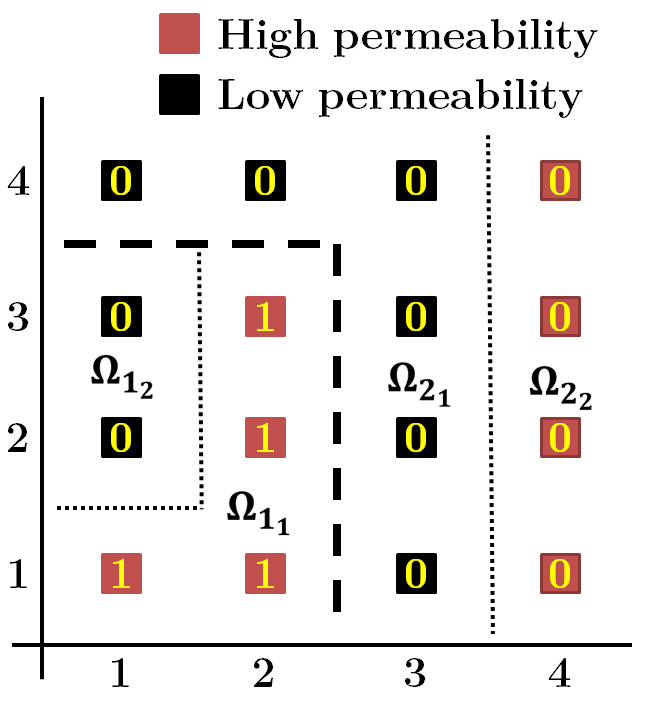}
						\label{fig:subdomain_levelset} }					
    \end{center}
    \caption{Subdomain (a), levelset (b) and subdomain-levelset deflation (c).}
    \label{fig:subdomain-levelset-deflation}
\end{figure}

If the discontinuities in the computational domain exhibit a complex geometry, subdomain-levelset deflation can be used to guarantee a good approximation of the eigenvectors corresponding to extreme eigenvalues \cite{Tang07b}. Whereas subdomain deflation does not take jumps into account, subdomain-levelset deflation identifies different regions in the domain with similar properties. A simple example is given in Figure \ref{fig:subdomain-levelset-deflation}.
The grid is $4 \times 4$ and nodes are shown as squares. In each case, the values shown on the nodes correspond to the values in the first deflation vector. In the middle and right figure, the border between the red nodes (high permeability) and black nodes (low permeability) exemplifies a sharp contrast in the PDE coefficient. The figures shows the following:
\begin{itemize}
	\item In the left figure, subdomain deflation is used. The dashed line divides the domain into the four subdomains $\Omega_1, \Omega_2, \Omega_3$ and $\Omega_4$. Each subdomain corresponds to a unique deflation vector.
	\item In the middle figure, levelset deflation is used. This time, the dashed line coincides with the contrast in the PDE coefficient. As a result, we get two domains $\Omega_1$ and $\Omega_2$.
	\item In the right figure, subdomain-levelset deflation is used. The subdomain division is determined using certain criteria, which in this example leads to the division (dashed line) between $\Omega_1$ and $\Omega_2$. Within each subdomain, levelset deflation (dotted line) uses the jump between the high permeability and low permeability nodes to obtain the subdomains $\Omega_{1_1}, \Omega_{1_2}, \Omega_{2_1}$ and $\Omega_{2_2}$
\end{itemize}

The pseudocode for deflated GMRES using the subdomain-levelset method is mostly identical to the harmonic Ritz deflation code given in Algorithm \ref{alg:gmres-defl-HR}, with two differences. First, instead of assigning $P_1 = P_2 = I$ and $x^* = 0$ in the first line, the physics-based deflation vectors are constructed (manually or automatically) a priori and we compute \eqref{eq:deflops} and $x^* = ZE^{-1}Z^Tb$ before the iteration starts. Because the initial permeability is fixed, $P_1$ and $P_2$ remain constant throughout the simulation. Second, as a result from the a priori construction of the deflation vectors, lines \ref{alg:gmres-defl-dyn:line:ifstart} - \ref{alg:gmres-defl-dyn:line:ifend} in Algorithm \ref{alg:gmres-defl-HR} can be omitted.

In this paper, only the serial implementation of the subdomain-levelset deflation method is used. We highlight the fact that the algorithm is particularly suitable for a parallel implementation. When the subdomain division is determined by the parallel partitioning, as is often the case in commercial reservoir simulation, each parallel subdomain corresponds to a processor. We can apply the levelset deflation method to each subdomain, and append the deflation vectors with zeros for all cells outside the parallel subdomain. Not only can the setup phase be executed nearly completely in parallel, the resulting set of (sparse) vectors is also particularly suitable for a parallel implementation of the deflation operators $P_1$ and $P_2$. For a fixed computational domain, the parallel subdomains will become smaller as we increase the number of processors. Hence, the levelset method can be applied to a smaller domain, which reduces the cost of the setup phase. In addition, more details might be captured by the (fixed number of) deflation vectors.

Whereas previous work in this area, such as \cite{Frank01,Tang05}, has mainly been directed towards the parallel subdomain deflation method, we have implemented a parallel subdomain-levelset deflation method. The algorithm to construct the deflation vectors automatically, and the parallel implementation of the deflation operators, will be discussed in an upcoming publication. In this work, we will construct the deflation vectors manually, using the guidelines of the levelset deflation method.

\section{Numerical experiments}
Harmonic Ritz and subdomain-levelset deflation are implemented both in Matlab and a full simulation in C++. The former runs only a single pressure solve, while the latter includes a non-linear solver (Newton-Raphson) and a CPR-preconditioner. The pressure solve, in this case, is the first stage of the CPR preconditioner. By default, the pressure solve is preconditioned by AMG. The second stage of CPR uses ILU. For a detailed discussion of the implementation, we refer to \cite{vanderLinden13}. 

The experiments are run on three cases, varying in size and complexity. In the upcoming subsections, each case will be discussed in terms of dimensions, initial conditions for the permeability and wells. A few simulation times are picked as reporting times $T$ (in days).

\begin{figure}[H]
\centering
\includegraphics[width=5cm]{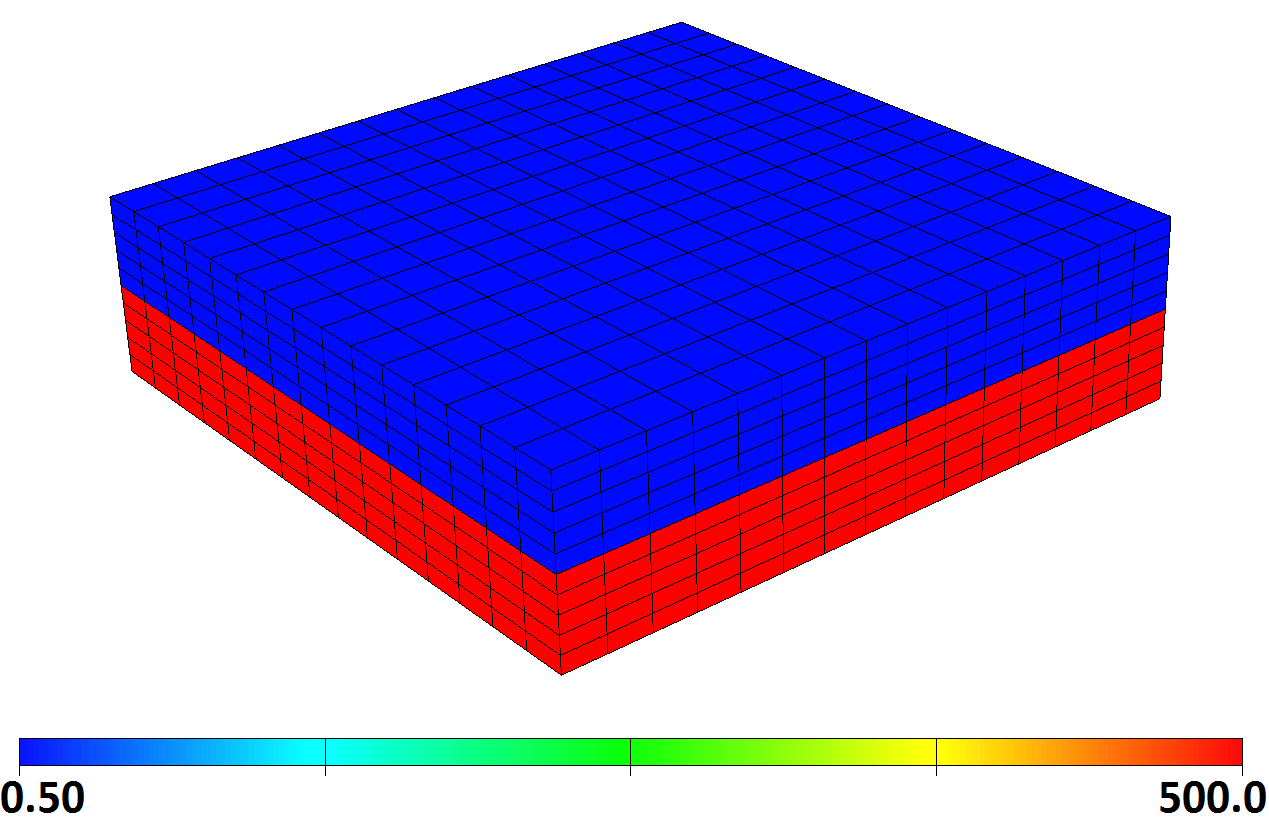}
\caption{Permeability field of the BO case.}
\label{fig:BO-permXYZ}
\end{figure}

\textbf{Model A: Black Oil Model}
The Black Oil (BO) model is the first case under consideration.  The size of the BO case is $15 \times 15 \times 10$. The initial permeability (equal in $x,y,z$ directions) is shown in Figure \ref{fig:BO-permXYZ}. 
A relatively small permeability jump is present between the fifth (red) and sixth (blue) horizontal layer of cells. The model contains nine wells: seven producers and two injectors. The injectors inject water at a fixed rate throughout the simulation. 

\begin{figure}[H]
\centering
\includegraphics[width=15cm]{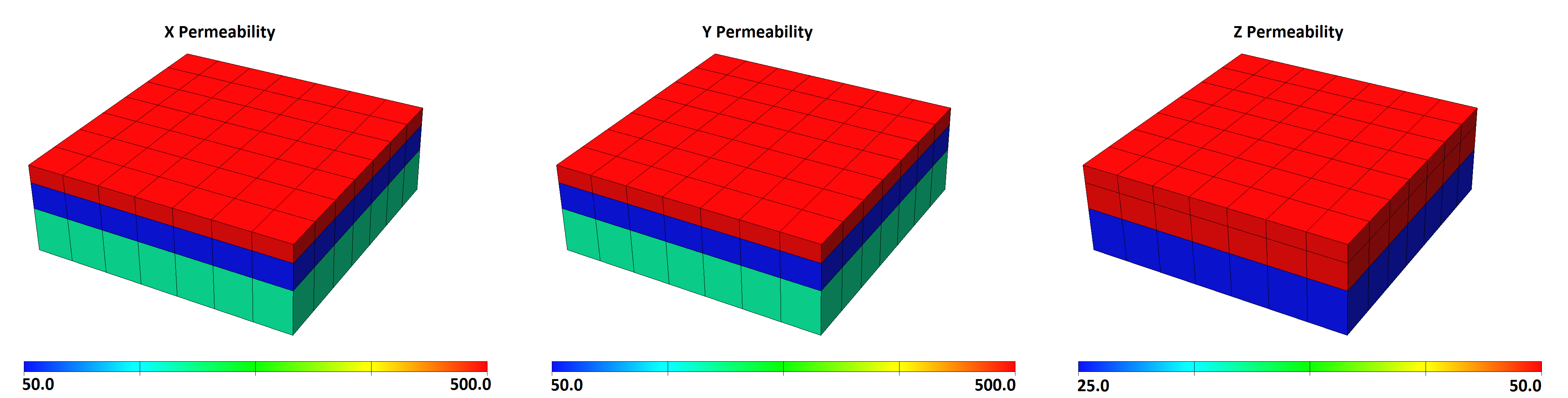}
\caption{Permeability in $x$, $y$ and $z$-direction.}
\label{fig:SPE5-XYZperm}
\end{figure}

\textbf{Model B: SPE fifth comparative solution project}
The fifth comparative solution project of the Society of Petroleum Engineers (SPE5) is part of a series of comparative solution problems designed to compare reservoir simulators from different companies, research institutes and consultants in the petroleum industry \cite{Killough87}. SPE5 focuses on the simulation of the (miscible) flooding of a reservoir. The dimensions of the SPE5 case are $7 \times 7 \times 3$. The default initial permeability in the $x$, $y$ and $z$ direction is shown in Figure \ref{fig:SPE5-XYZperm}.
In vertical direction, the permeability jump between the second and third layer is small. Within each horizontal layer in $x$ and $y$ direction, the permeability is constant. A modified version of the SPE5 case was already introduced in Figure \ref{fig:sandwich}. One injector and one producer are placed in opposite corners of the reservoir. Starting at $T = 0$ days, the injector pumps water in the reservoir, which pushes the oil towards the producer. 

\textbf{Model C: Steam Assisted Gravity Drainage}
A technique called Steam Assisted Gravity Drainage (SAGD) uses steam injection to create a steam chamber around the producers. The reservoir is heated to make the oil less viscous, after which water is injected into the reservoir. The water evaporates to become steam, creating the steam chamber. After expanding the steam chamber upwards, gravity causes the heavy oil to flow down to the production wells. Due to the steam injection and temperature effects, large pressure gradients occur around the injector and the producer. The gradients, in turn, produce strong non-linearities in the solution of the SAGD case.

\begin{figure}[H]
\centering
\includegraphics[width=15cm]{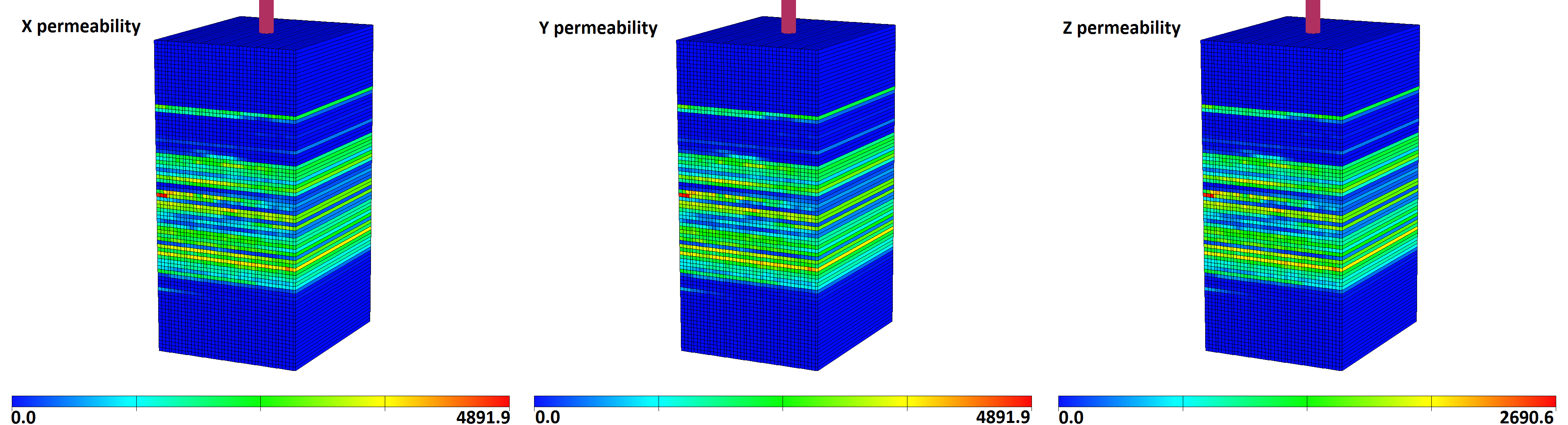}
\caption{Permeability in $x$, $y$ and $z$ directions.}
\label{fig:SAGD-perm}
\end{figure}

The SAGD-SMALL case has dimensions $41 \times 1 \times 85$. The initial permeability for the SAGD case is the same in $x$ and $y$ direction. For the $z$ direction, the pattern is the same but the permeability values are halved. Figure \ref{fig:SAGD-perm} shows the horizontal layer structure that is commonly found in petroleum reservoirs.
The top and bottom half of the reservoir have zero permeability. Near the center, we see several large permeability jumps of order $10^3$. Two wells, both capable of acting as a producer and an injector, are placed above each other near the bottom of the reservoir, allowing for the steam assisted gravity drainage. The sudden temperature and pressure changes render this case relatively difficult to solve. 

\textbf{Model D: SPE tenth comparative solution project}
In addition, we investigate the use of physics-based deflation in the large scale simulation. For this reason we
consider a well-known SPE-10 benchmark \cite{Christie01}. The fine grid is 180 x 220 x 85 for a total of 3366000 
cells with $\mathcal{O}(10^6)$ degrees of freedom
\begin{figure}[H]
\centering
\includegraphics[width=8cm]{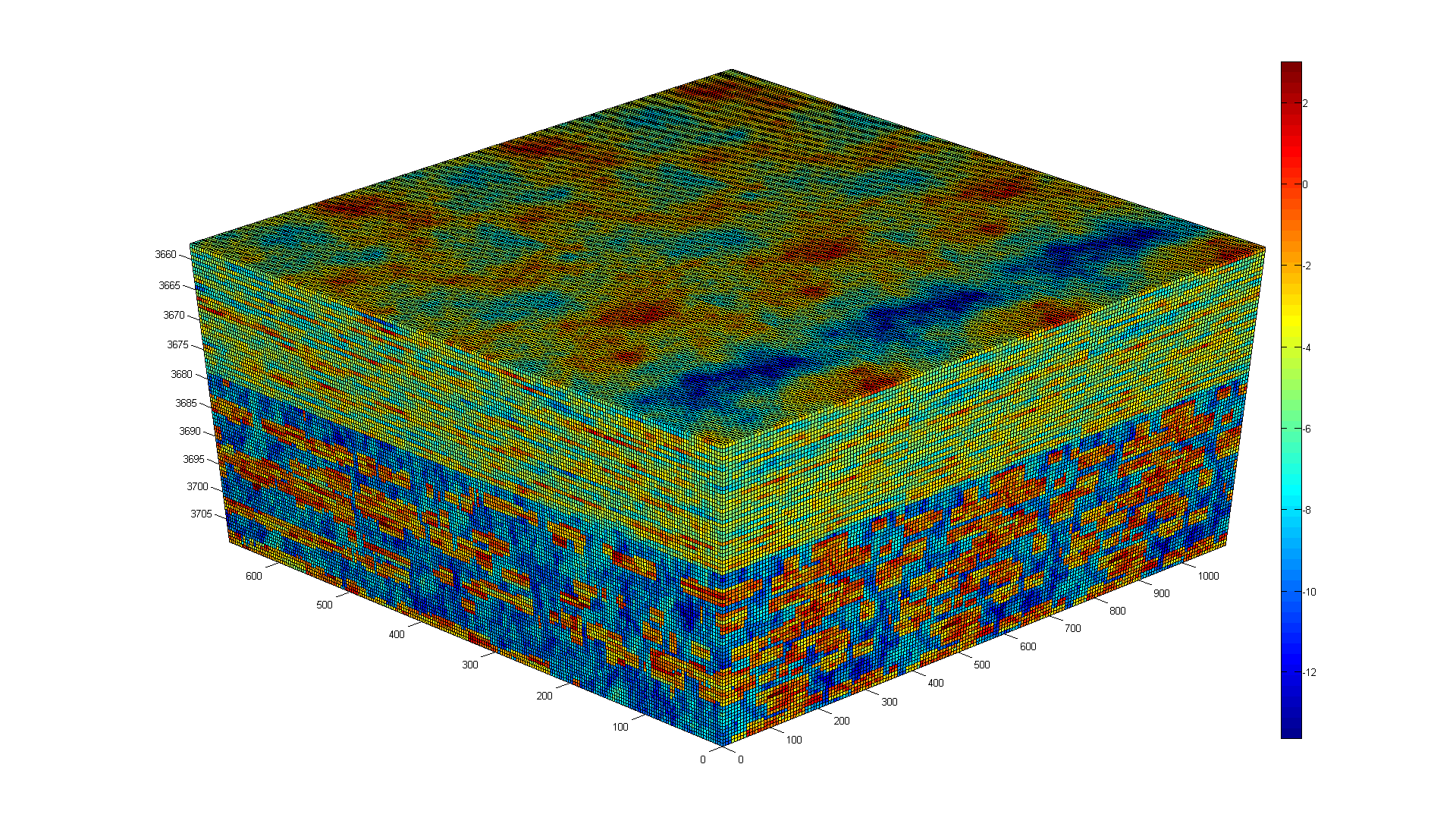}
\caption{X-direction permeability distribution in logarithmic scale and Darcy units for SPE10 case.}
\label{fig:SPE10}
\end{figure}

\subsection{Harmonic Ritz deflation}
\label{sec:hritzres}
In this section, we discuss the results obtained from applying harmonic Ritz deflation to the BO, SPE5 (with modified permeability) and SAGD-SMALL cases. The implementation in Algorithm \ref{alg:gmres-defl-HR} is used both in Matlab and the full simulation. 
\begin{rmk}
	Deflated GMRES using exact eigenvectors is denoted as DGMRES($m$,$d$), where $m$ is the cycle size and $d$ is the number of deflation vectors used. Harmonic Ritz deflation is denoted as RDGMRES($m$,$d$). Unless noted otherwise, we assume that the deflation vectors correspond to the smallest $d$ (approximated) eigenvalues.
\end{rmk}
In most of the experiments in this section, a simple Jacobi preconditioner is used, applied from the right. We use Jacobi because it is relatively cheap.

\subsubsection{Matlab simulation}
\begin{figure}[H]
\centering
\includegraphics[width=7cm]{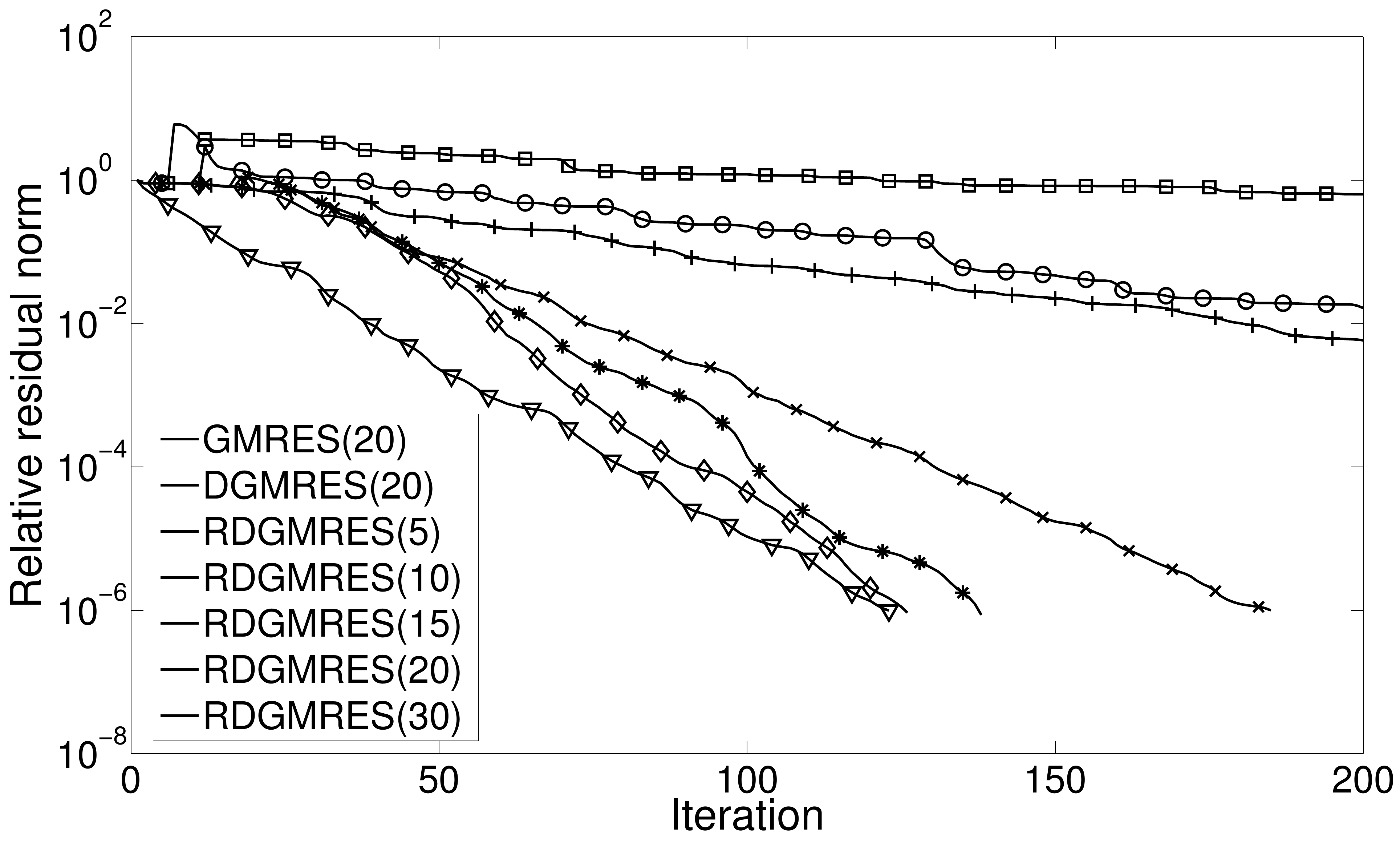}
\caption{DGMRES and RDGMRES convergence residuals for varying $m$.}
\label{fig:BO98-hritzdefl}
\end{figure}

\begin{figure}[H]
     \begin{center}
        \subfigure[] {
            \includegraphics[width=7cm]{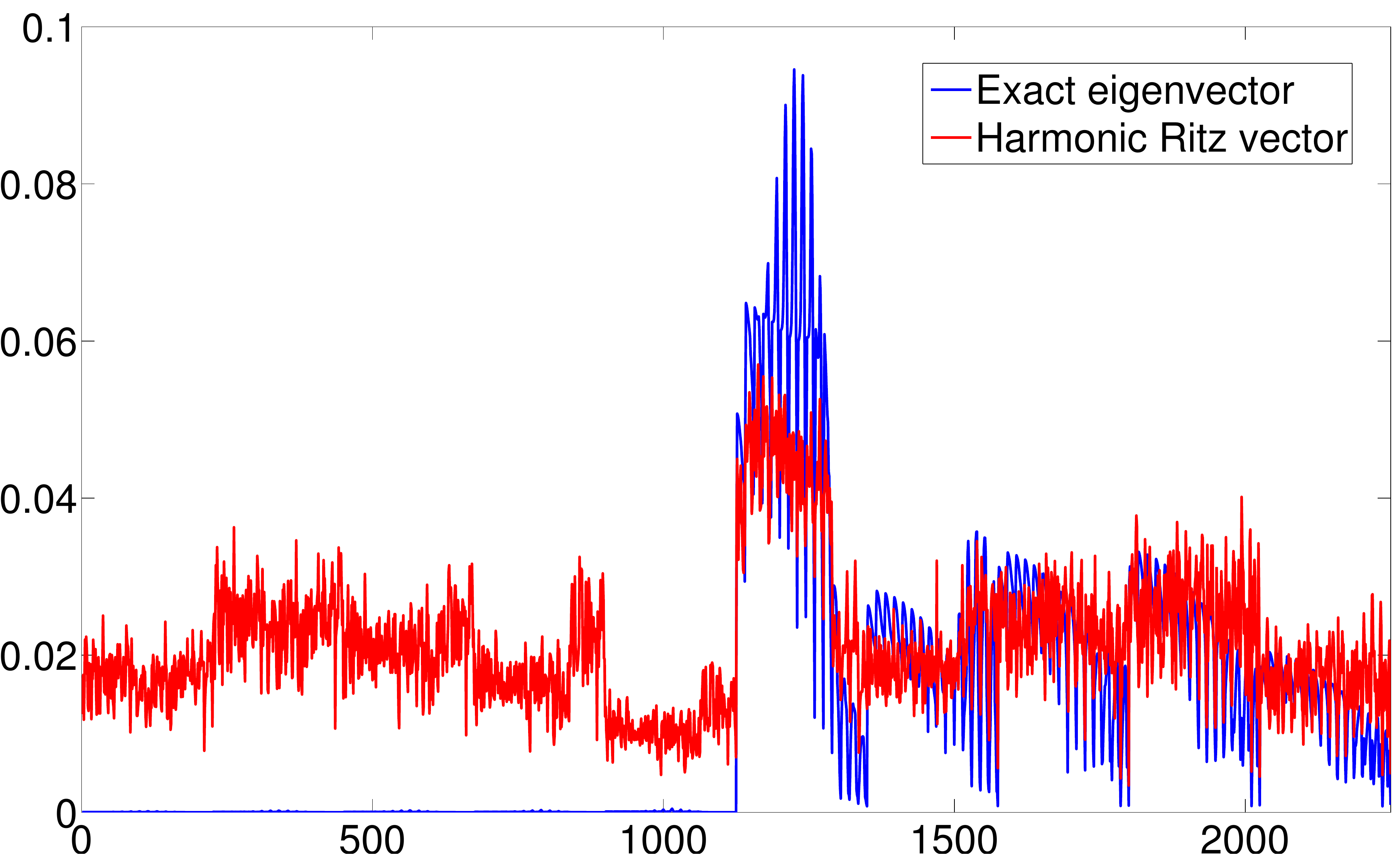}
						\label{fig:BO98_res5ritz} }
				\subfigure[] {
            \includegraphics[width=7cm]{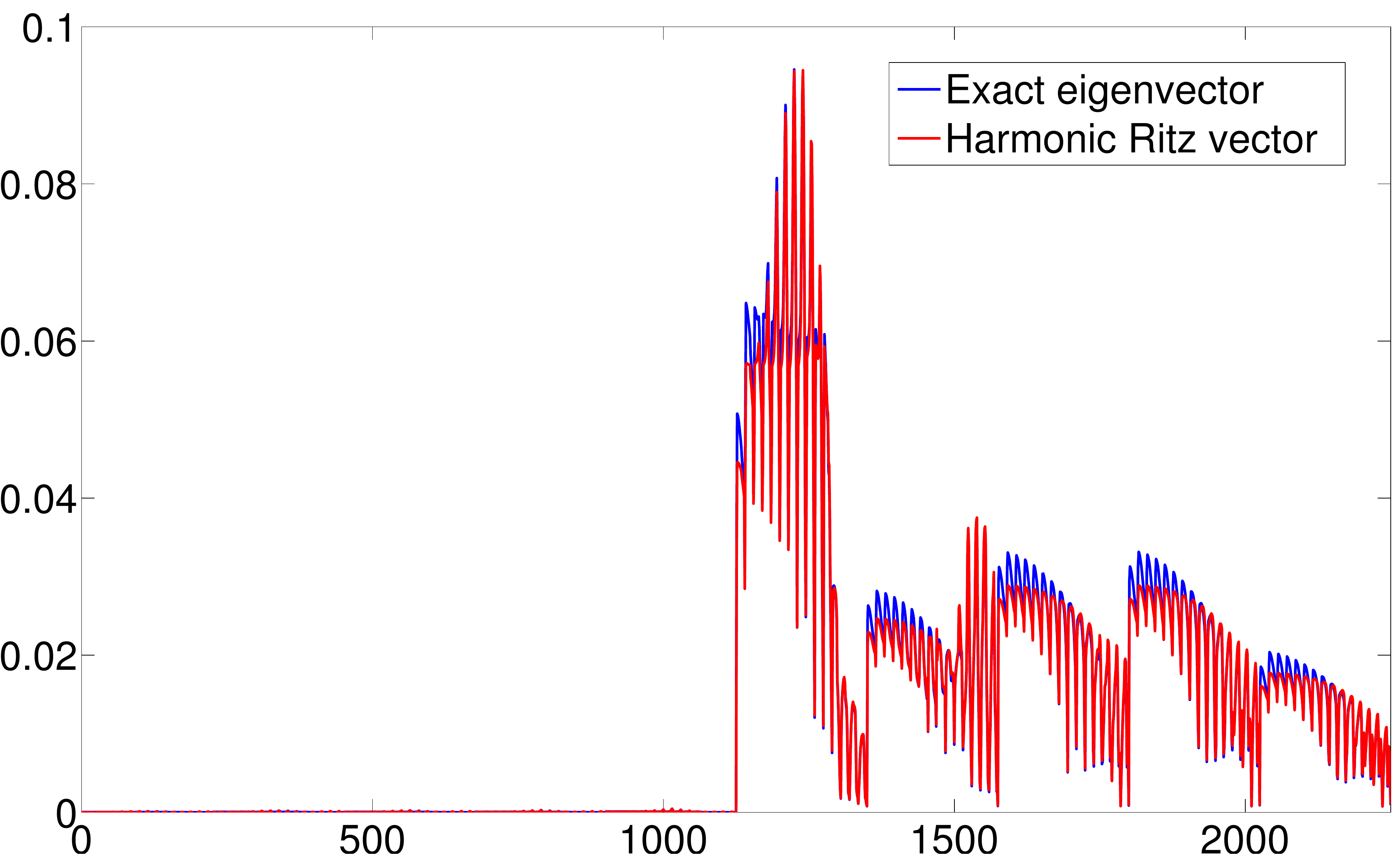}
						\label{fig:BO98_res30ritz} }						
    \end{center}
    \caption{The harmonic Ritz vector and the exact eigenvector after (a) five and (b) thirty iterations.}
    \label{fig:BO98-res530ritz}
\end{figure}

In Figure \ref{fig:BO98-hritzdefl}, we plot the linear solve of the BO case using GMRES(20), DGMRES(20,1) using one exact eigenvector and RDGMRES($m$,1) using one harmonic Ritz (HR) vector. The cycle size $m$ is varied to demonstrate the impact on convergence.
The fastest convergence is achieved by DGMRES(20,1). RDGMRES(30,1) has higher residual norms in the first 100 iterations, but reaches the tolerance of $10^{-6}$ at nearly the same iteration count. For lower values of $m$, the convergence is slower, although still faster than GMRES without deflation. For $m = 10$ and $m = 5$, however, the harmonic Ritz vector is no longer a sufficiently accurate approximation of the true eigenvector. As a result, RDGMRES is slower than GMRES. To illustrate this phenomenon, the exact eigenvector and harmonic Ritz eigenvector approximation after five and thirty iterations are plotted in Figure \ref{fig:BO98-res530ritz}.
Clearly, the harmonic Ritz deflation vector is not a good approximation of the true eigenvector after five iterations, resulting in poor convergence of RDGMRES. Equation \eqref{eq:deflres} does not hold, because the deflation subspace is not (near) invariant (see Section \ref{sec:gmresconv}). After thirty iterations, on the other hand, the harmonic Ritz vector and the eigenvector nearly overlap and deflation significantly improves the convergence of GMRES. For intermediate choices of $m$, the approximation will be more/less accurate, as demonstrated in the convergence history in Figure \ref{fig:BO98-hritzdefl}. 

After repeating the above experiment for the SPE5 case with modified permeability, we find that the size of the permeability jump imposes a requirement on the cycle size. If we run the simulation with permeability jump $\sigma = 10^8$, for example, it turns out that $m = 20$ is not sufficient to obtain sufficiently accurate harmonic Ritz eigenvector approximations. As a result, RDGMRES will not converge. For this particular case, $m \geq 40$ is required for deflation to be effective. We believe that this correlation is caused by the extreme eigenvalues in the SPE5 case, which become more isolated and extreme as we increase $\sigma$. The more extreme the eigenvalues become, the longer it takes for the smallest Ritz values to converge to the extreme eigenvalues. Consequently, more iterations are required in a cycle to find decent approximation of the eigenvectors.

\subsubsection{Full simulation}
Having demonstrated the potential of the harmonic Ritz deflation method in Matlab, we continue with the results for full simulation. Each simulation in this section uses $A_p = A$ (see Algorithm \ref{alg:gmres-defl-HR} and preconditioning applied from the right). All simulations are run in serial. In the previous section, we showed that harmonic Ritz deflation can effectively eliminate the harmful eigenvalues and thereby improve convergence. As will become clear in this section, however, the performance in the full simulation is not as good and the overhead of deflation is too high to compete with the default CPR preconditioning scheme.

The following variables are used to compare the results:
\begin{itemize}
	\item \textbf{Non-linears}. The amount of non-linear iterations in the simulation.
	\item \textbf{Fails}. The amount of failed non-linear iterations.
	\item \textbf{Outer linears}. The amount of iterations used to solve the linear systems generated by the non-linear iterations.
	\item \textbf{Inner linears}. The amount of iterations used to solve the pressure systems generated by the second stage of the CPR preconditioning in the outer GMRES loop.
	\item \textbf{CPU time}. The overall CPU time (in seconds) of the linear solve.
\end{itemize}
We refer to the linear solve of the full linear system and the pressure solve as the outer linear solve and the inner linear solve, respectively. For each outer linear iteration, there is at least one inner linear iteration. Note that the number of non-linears is equal to the number of linear systems that need to be solved in the outer linear solve. Similarly, the number of outer linears is equal to the number or pressure systems that need to be solved in the inner linear solve. We compare the overall CPU time of the linear solve to analyze the overhead of deflation. 

For the pressure solve in commercial simulations, often only a single inner linear iteration is performed, using an AMG preconditioner. Harmonic Ritz deflation cannot be used with a single-iteration pressure solve, hence the following settings are used for our experiments in this section.
\begin{table}[H]
	\centering
	\begin{tabular}{@{}ll@{}}
		\midrule[1.2pt]
		 Setting &  Default value \\
		\midrule
		Pressure solve tolerance & $10^{-6}$ \\
		Cycle size & $30$ \\
		Minimum number of iterations & $30$ \\
		Maximum number of iterations & $60$ \\
		Preconditioner & Jacobi \\
		\midrule[1.2pt]
	\end{tabular}
	\caption{Settings for the results in Table \ref{tab:deflIX}.}
	\label{tab:exp1IX}
\end{table}
In order to use harmonic Ritz deflation, a restart is required. Therefore, we set the minimum number of iterations equal to the cycle size. By using a weaker (but cheap) Jacobi preconditioner, we highlight the advantage of deflation. The cycle size is chosen sufficiently large to guarantee accurate eigenvector approximations. For moderate heterogeneity in the permeability, our experience is that $m$ should be at least $20$. 

\begin{figure}[H]
\centering
\includegraphics[width=8cm]{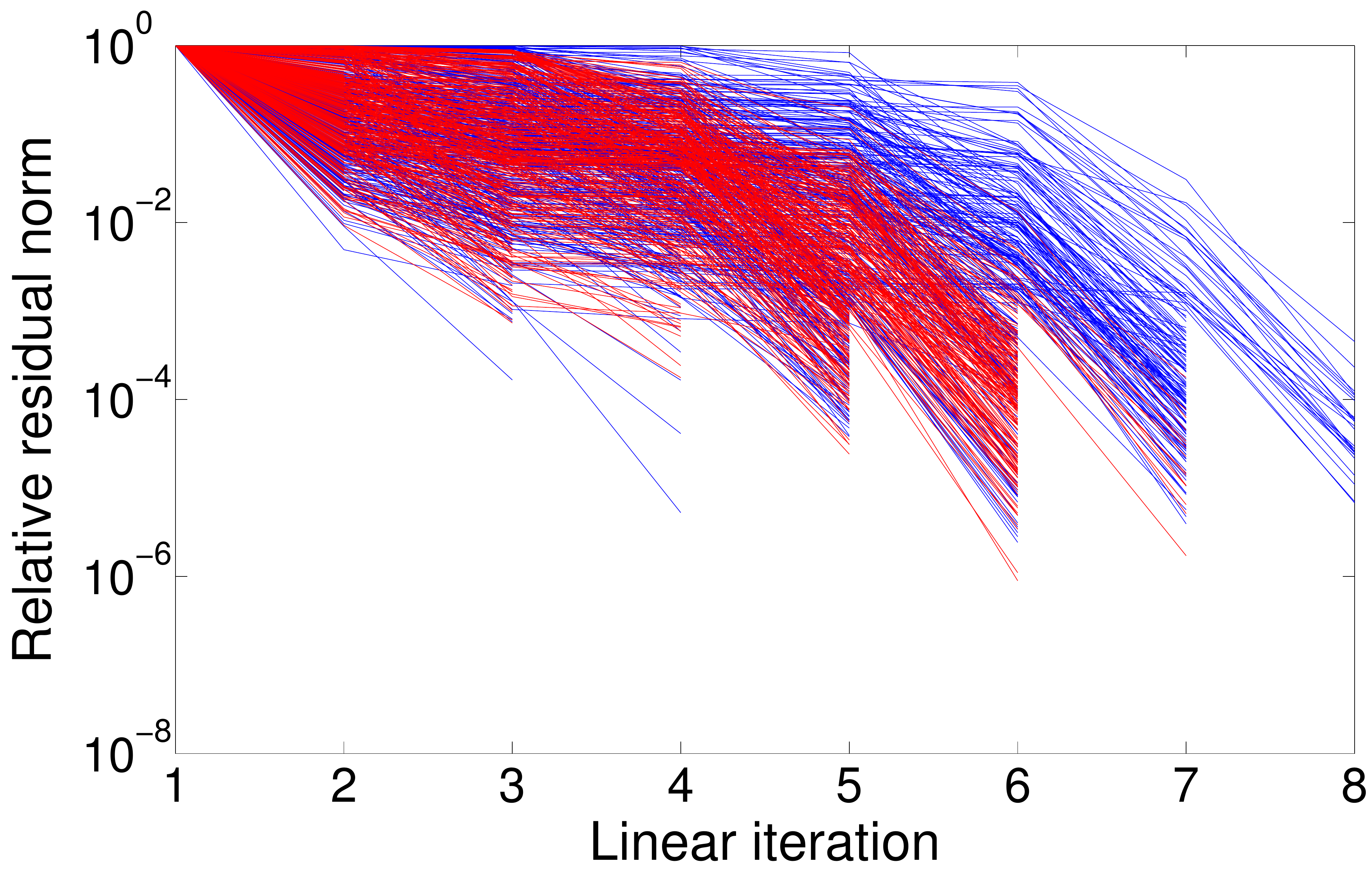}
\caption{Inner linear convergence of GMRES(30) with Jacobi (in blue) and DGMRES(30,3) with Jacobi (in red).}
\label{fig:SPE5-IX-3ritzdeflfrozen}
\end{figure}

Table \ref{tab:deflIX} shows a comparison of GMRES(30) and RDGMRES(30,3) in the two left columns, using the settings from Table \ref{tab:exp1IX} and the SPE5 case with modified permeability ($\sigma = 10^8$). The convergence of all inner linear iterations in the GMRES(30) and RDGMRES(30,3) pressure solve (both with Jacobi preconditioner) is plotted in Figure \ref{fig:SPE5-IX-3ritzdeflfrozen}.
\begin{table}[H]
	\centering
	\begin{tabular}{@{}llll@{}}
		\midrule[1.2pt]
		 & GMRES(30) & RDGMRES(30,3) & GMRES(30) \\
		\midrule
			Preconditioner		   & Jacobi & Jacobi & AMG \\
			Non-linears            & $527$ & $430$ & $395$\\
			Fails                  & $0$ & $0$ & $0$\\
			Outer linears          & $2,057$ & $1,606$ & $1,099$\\
			Inner linears          & $123,420$ & $96,360$ & $6,659$\\
			CPU time & $6.58$ & $5.68$ & $1.22$\\
		\midrule[1.2pt]
	\end{tabular}
	\caption{Comparison of GMRES(30) (with Jacobi preconditioner), RDGMRES(30,3) (with Jacobi preconditioner) and GMRES(30) (with AMG preconditioner) for the modified SPE5 case.}
	\label{tab:deflIX}
\end{table}
Firstly, note that $2,057 \times 60 = 123,420$ and $1,606 \times 60 = 96,360$, which implies that neither GMRES(30) (with Jacobi) nor RDGMRES(30,3)  reaches the tolerance of $10^{-6}$ in any of the pressure solves. The number of outer linears is reduced by approximately $22\%$ because the residual in each pressure solve is (slightly) better for RDGMRES compared to GMRES. As a result, the outer linear solve will converge slightly faster. From the CPU time we conclude that the overhead of deflation is small enough to render the method efficient. Although the speedup is not as significant as was shown in the Matlab experiments in Figure \ref{fig:BO98-hritzdefl}, deflation in the diagonally scaled pressure system is still efficient. From Figure \ref{fig:SPE5-IX-3ritzdeflfrozen}, it is clear that the convergence has improved by using deflation.

By inspecting the convergence history, we find that GMRES(30) with AMG reaches the tolerance in the pressure solve after (on average) $6,659/1,099 \approx 6$ inner linear iterations. RDGMRES with Jacobi requires over $14$ times as many total inner iterations, and is $4.7$ times slower in terms of CPU time than GMRES with AMG. Moreover, in the default settings the pressure solve uses only one AMG-preconditioned GMRES iteration. Using more than a few iterations becomes too costly when using very large cases with millions of nodes. Our experiments show that the residual in this pressure solve setup is often sufficient to guarantee good convergence of the outer linear solve. When only one iteration of GMRES is used, the CPU time of the linear solve will be a fraction of the $1.22$ seconds above.

For this experiment, it is not possible to use RDGMRES with an AMG preconditioner. As noted before, at least 20 iterations are required to obtain sufficiently accurate approximations of the eigenvectors. AMG-preconditioned GMRES, however, reaches machine-precision residual norms within 10 to 15 iterations. The physics-based deflation method in the next section is more suitable to be combined with AMG, because the physics-based deflation vectors are applied from the start of the simulation. Our experiments for this deflation method show that adding deflation to AMG-preconditioned GMRES does not always improve convergence. We believe that, in some cases, AMG is already capable of tackling the harmful eigenvalues. In \cite{Trottenberg01}, the occurrence of extreme eigenvalues is linked to `algebraically smooth' error nodes. AMG is very efficient in reducing these error characteristics in general, which explains why adding deflation does not always improve convergence.

\subsubsection{Optimal RDGMRES settings}
In an attempt to reduce the overhead of harmonic Ritz deflation, we have tried reducing the cycle size and minimum number of iterations to 20, while keeping the maximum number of iterations at 60. Furthermore, we know from our simulations that the water injection in the SAGD-SMALL case causes issues in simulation. Shortly after the water injection starts, the non-linear time step size is reduced several times before the simulation can continue. Because isolated eigenvalues may be responsible, we tried switching on deflation when the water injections starts. To reduce the computational costs, we did not use deflation for all previous time steps. Lastly, we varied the number of deflation vectors.

None of the attempts produced satisfactory results. Using deflation on the more difficult SAGD-SMALL case leads to a 17\% inner linear iteration reduction, compared to 5\% in the SPE5 case. The overhead increases, however, due to the increased size of the system matrix. Switching deflation off up to the time of the water injection decreases the overhead, yet the gain is offset by the increased number of linear iterations. In general, the amount of inner linears decreases as we increase the number of deflation vectors. The time gain from the decreased number of inner linears does not outweigh the overhead of deflation, however, as the overall CPU time is higher if we use more deflation vectors.

In conclusion, none of the approaches lead to improvements. We formulate the following two (related) objectives for an improved deflation method:
\begin{itemize}
	\item The overhead of harmonic Ritz deflation offsets the time gain from the improved convergence. Therefore, we require a cheaper deflation method.
	\item The harmonic Ritz deflation vectors do not speed up convergence enough to offset the overhead. Therefore, we require more effective deflation vectors.
\end{itemize}
This conclusion led to the development of a physics-based deflation method. The advantages of this method have been briefly highlighted in Section \ref{sec:domdefl}, but will be repeated in the next section along with the results.

\subsection{Physics-based deflation}
\label{sec:pbd}
In physics-based deflation, we approximate the eigenvectors using the underlying physics. As discussed in Section \ref{sec:domdefl}, the span of the eigenvectors corresponding to the extreme eigenvalues in linear systems with strong heterogeneity can be approximated by the span of a set of physics-based deflation vectors. In our applications, the permeability is generally responsible for the largest jumps in the coefficients of the reservoir equations. The deflation vectors can be constructed manually, or computed automatically using a subdomain-levelset algorithm. The main advantages of physics-based deflation are:
\begin{itemize}
	\item The deflation vectors are defined a priori. Therefore, harmful eigenvalues can be eliminated from the spectrum from the start of the linear solve.
	\item The deflation vectors are computed only once and can be reused throughout the simulation.
	\item The method is relatively cheap and scalable. The setup can be executed nearly completely in parallel and the resulting vectors are sparse. Furthermore, because each deflation vector is zero outside the (parallel) subdomain, deflation can be implemented in parallel with limited communication.
	\item If regions of constant permeability are contained and separated by large, well-defined jumps, then constructing a set of efficient deflation vectors automatically is relatively easy. 
\end{itemize}
\subsubsection{Matlab simulation}
\begin{figure}[H]
     \begin{center}
        \subfigure[] {
            \includegraphics[width=4.5cm]{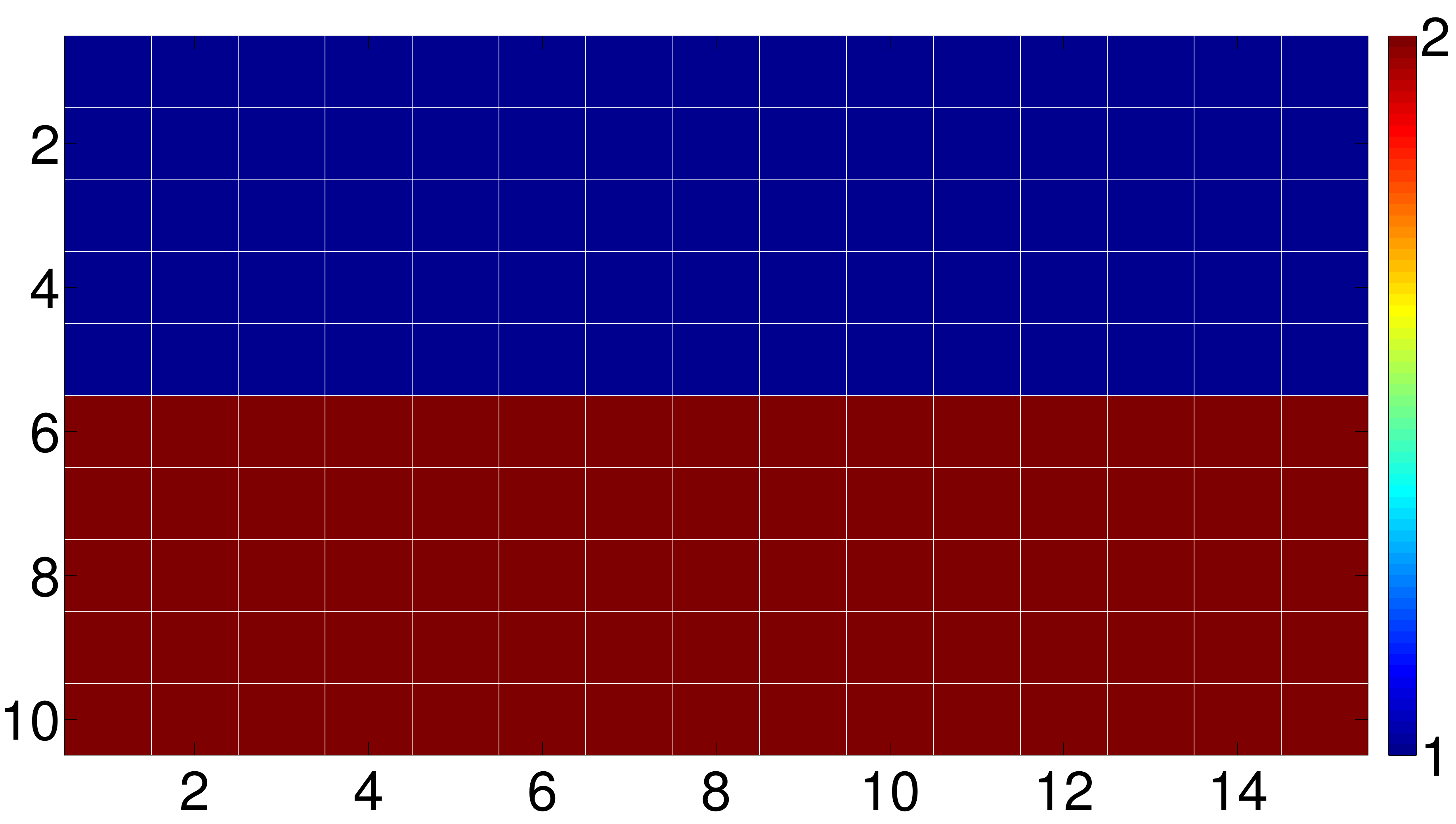}
						\label{fig:Mandefl_BO98_vecs} }
		\subfigure[] {
            \includegraphics[width=4.5cm]{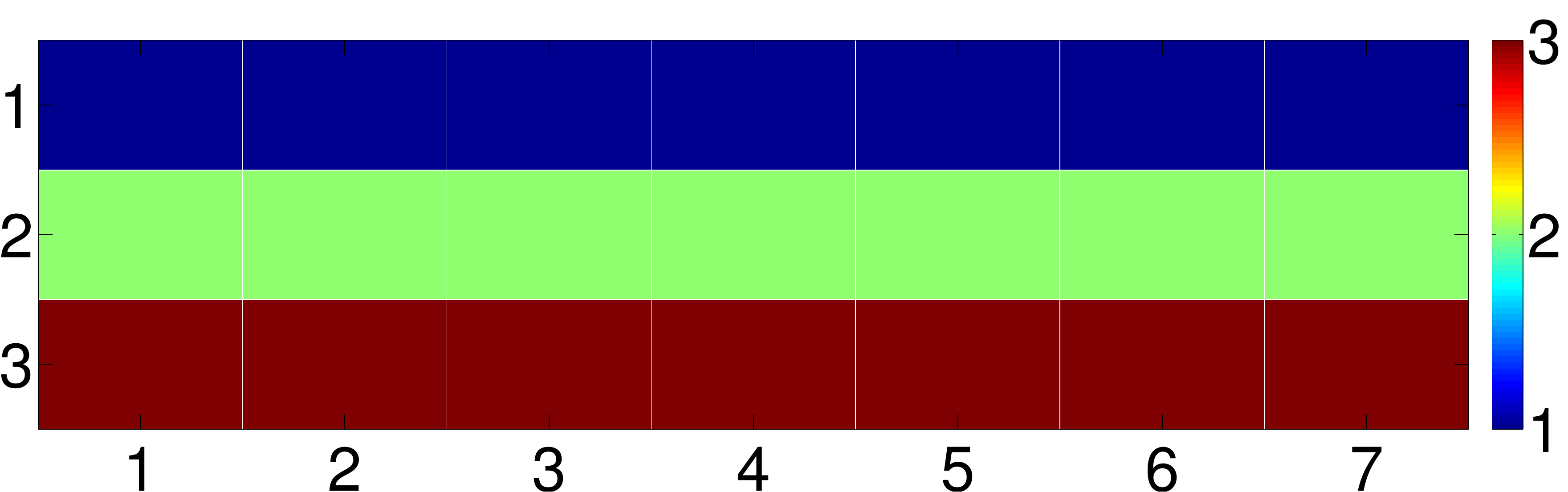}
						\label{fig:Mandefl_SPE5_vecs} }	
		\subfigure[] {
			\includegraphics[width=4.5cm]{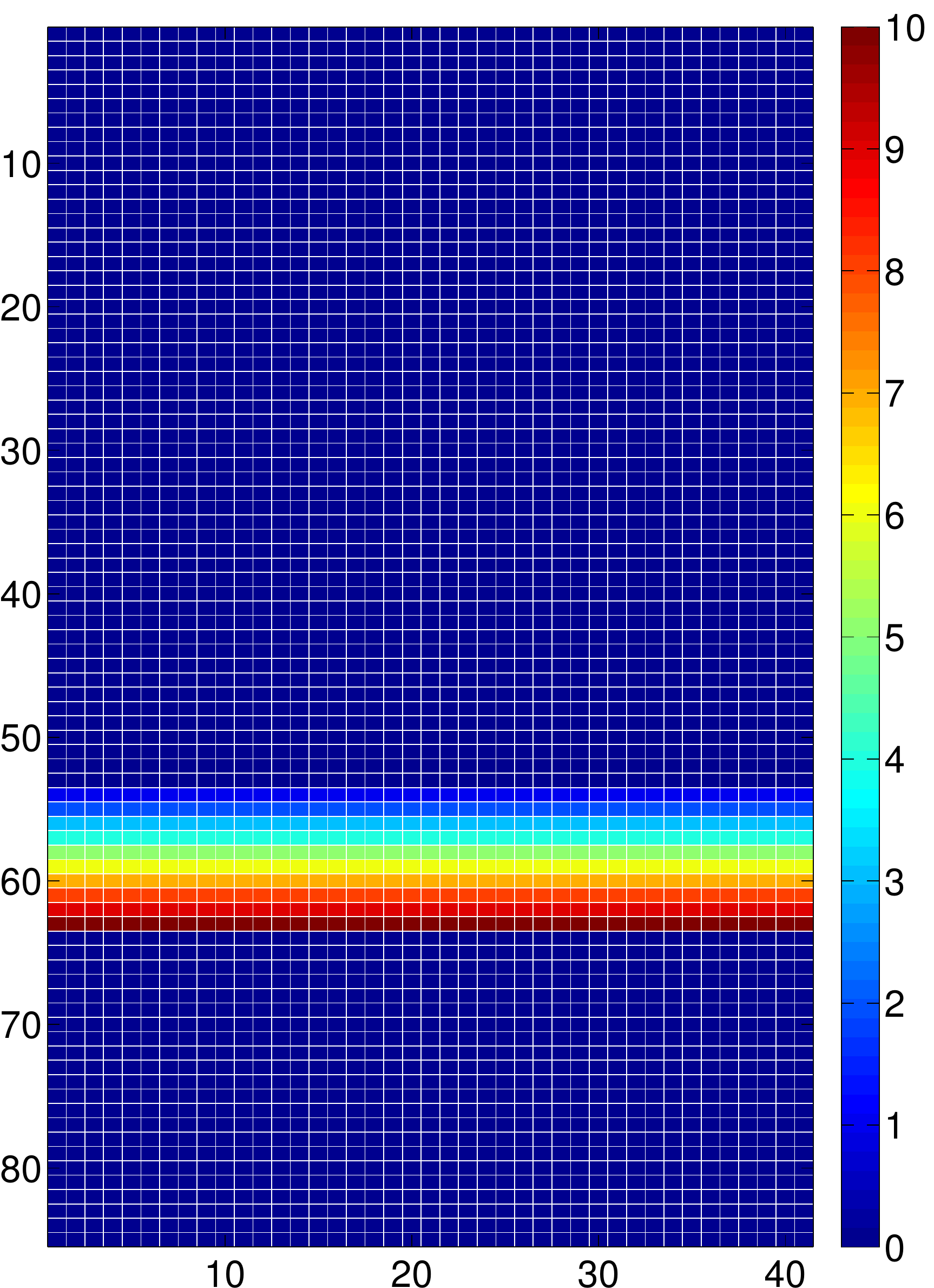}
						\label{fig:Mandefl_SAGD_vecs} }						
    \end{center}
    \caption{Manually constructed deflation vectors for (a) BO, (b) SPE5 and (c) SAGD-SMALL.}
    \label{fig:Mandefl-vecs}
\end{figure}

In the first set of numerical experiments with physics-based deflation, we manually assign deflation vectors. In particular,
\begin{itemize}
	\item The BO case has two horizontal layers of constant permeability, as shown in Figure \ref{fig:BO-permXYZ}. To capture the jump, we manually construct two deflation vectors using Equation \ref{eq:domdefl}. In the first deflation vector, we assign value 1 to each cell in the top layer, and value 0 to each cell in the bottom layer. In the second deflation vectors, we assign value 0 to each cell in the top layer, and value 1 to each cell in the bottom layer. The boundary between the ones and zeros, in this case, coincides with the jump in permeability.
	\item The SPE5 case has three horizontal layers of permeability, as shown in Figure \ref{fig:SPE5-XYZperm}. The jump is captured using three deflation vectors, whose values are assigned in the same manner as above.
	\item The SAGD-SMALL case has a more complex permeability field. Figure \ref{fig:SAGD-perm} shows a number of horizontal layers of approximately constant permeability. Between these layers, rather large permeability jumps occur. Our experience is that deflation vectors near the injector/producer (where the flow occurs) are the most effective. Therefore, we assign ten deflation vectors to this region, where each vector represents a horizontal layer.
\end{itemize}
The proposed physics-based deflation vectors for the BO, SPE5 and SAGD-SMALL case are illustrated in Figure \ref{fig:Mandefl-vecs}(a), (b) and (c), respectively. In each figure, a front view is shown of the reservoir. Each color represents a deflation vector with ones in the nodes with that color, and zeros elsewhere.
We now apply the manually constructed deflation vectors to the pressure solve in Matlab. In each simulation, the pressure matrix is preconditioned with a Jacobi preconditioner, applied from the right. The residual tolerance is $10^{-6}$, and, unless note otherwise, we use $m = 20$. The convergence of physics-based deflation (PDGMRES) is compared to harmonic Ritz deflation (RDGMRES) and GMRES without deflation. We use 2 deflation vectors for the BO case, 3 deflation vectors for the SPE5 case and 10 deflation vectors for the SAGD-SMALL case. 

\begin{figure}[H]
     \begin{center}
        \subfigure[] {
            \includegraphics[width=7cm]{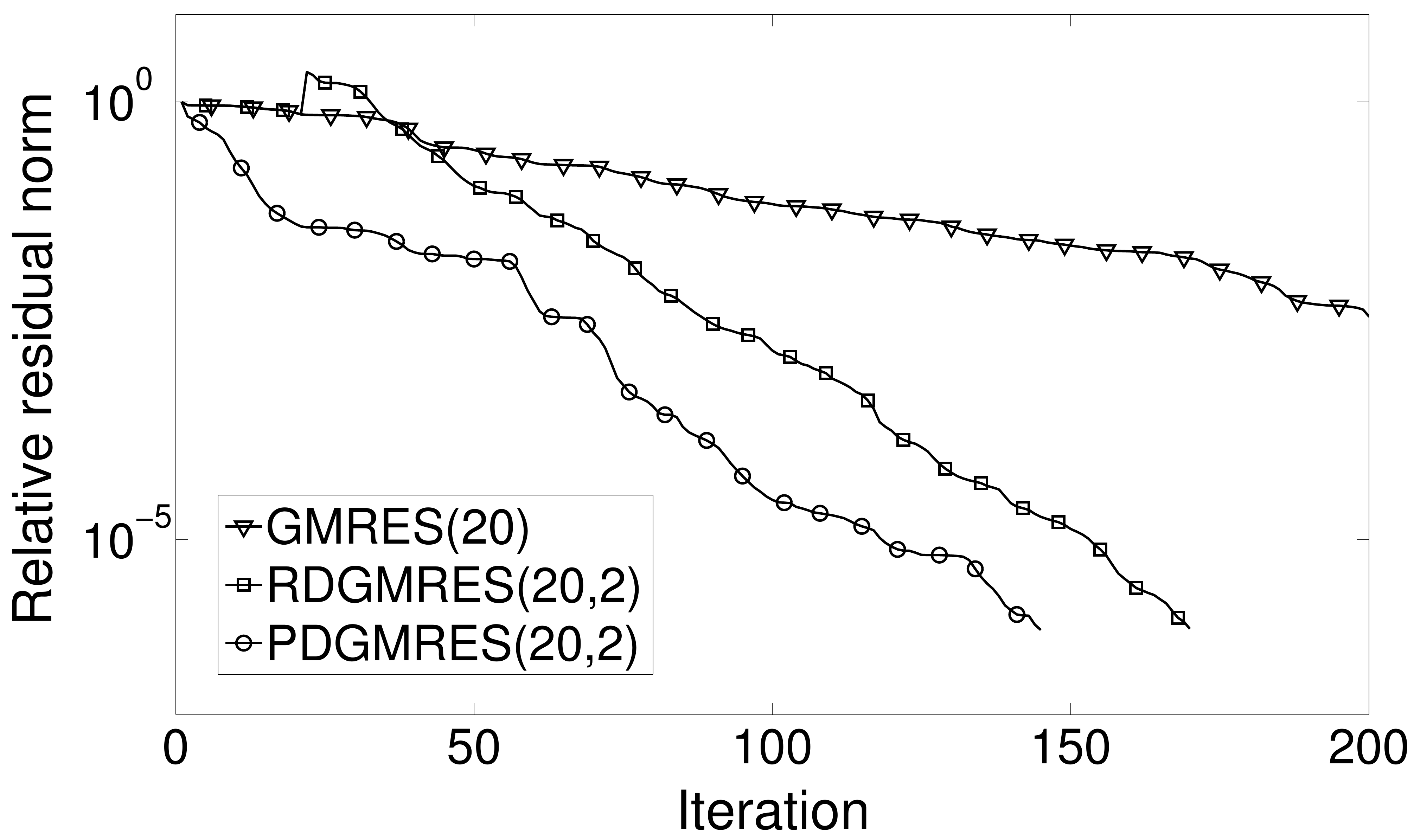}
						\label{fig:Mandefl_BO98_conv} }
				\subfigure[] {
            \includegraphics[width=7cm]{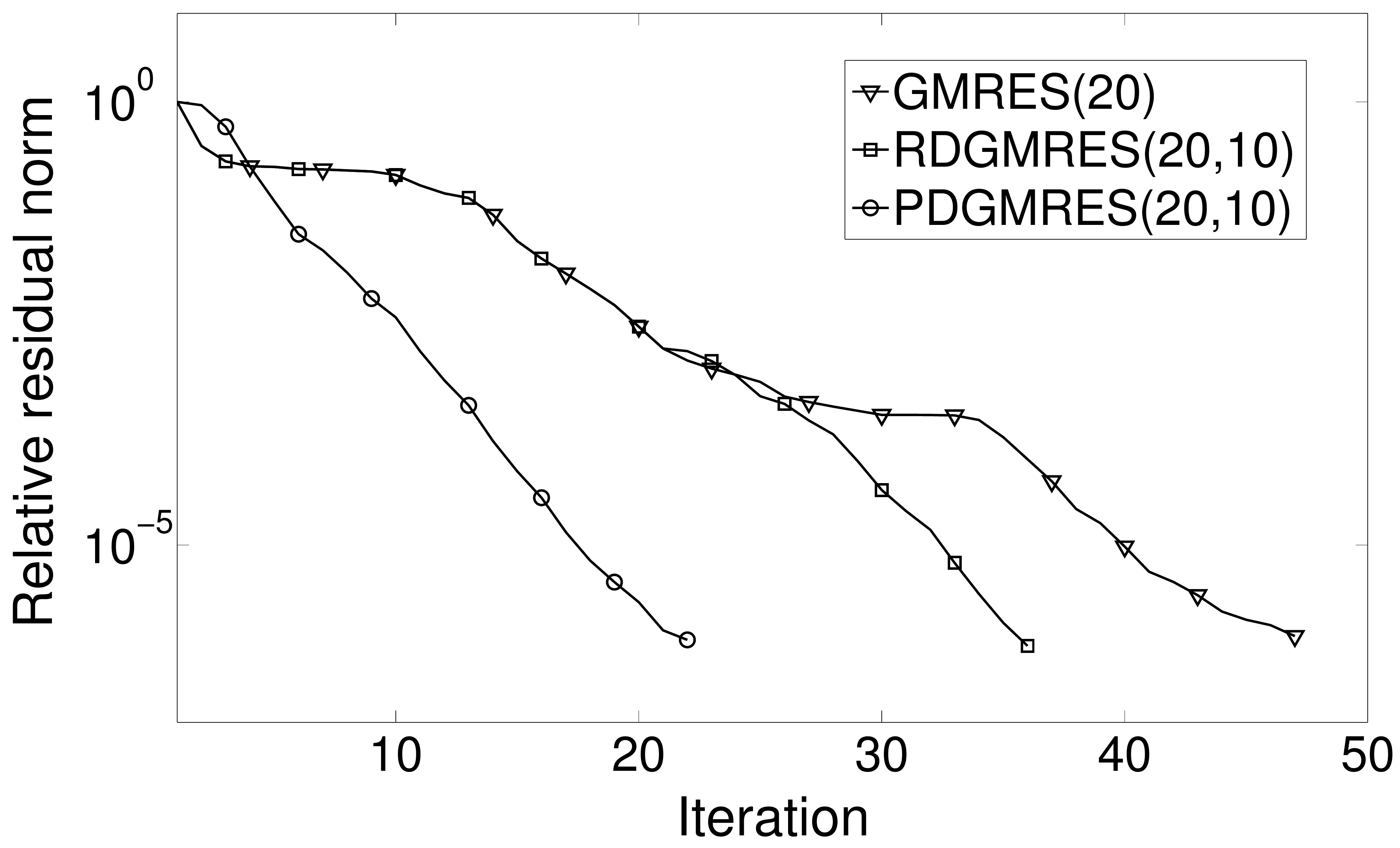}
						\label{fig:Mandefl_SAGD_conv} }						
    \end{center}
    \caption{Comparison of no deflation, harmonic Ritz deflation (RDGMRES) and (manual) physics-based deflation (PDGMRES) for (a) BO and (b) SAGD-SMALL.}
    \label{fig:Mandefl-BOSAGD}
\end{figure}

The convergence history of the BO case and the SAGD-SMALL case is shown in Figure \ref{fig:Mandefl-BOSAGD}.
In Figure \ref{fig:Mandefl-BOSAGD}(a), GMRES does not reach the tolerance within 200 iterations. The convergence of RDGMRES is better, but a residual increase occurs after the restart. The best residual convergence is attained by physics-based deflation. After 20 iterations, the convergence speed of PDGMRES is approximately equal to the convergence speed of RDGMRES, however physics-based deflation can be applied from the start of the simulation, which results in the best convergence. Similarly, PDGMRES in Figure \ref{fig:Mandefl-BOSAGD}(b) achieves the best performance. At the residual tolerance, the number of iterations compared to native GMRES is approximately halved. Again, the convergence speed of physics-based deflation after the restart is equal to the convergence speed of harmonic Ritz deflation, but the former method has the advantage of deflating the extreme eigenvalues from the start. Lastly, we note that the temporary residual increase for RDGMRES in Figure \ref{fig:Mandefl-BOSAGD}(a) can be related back to the discussion in Remark \ref{rmk:AM}.

\begin{figure}[H]
\centering
\includegraphics[width=8cm]{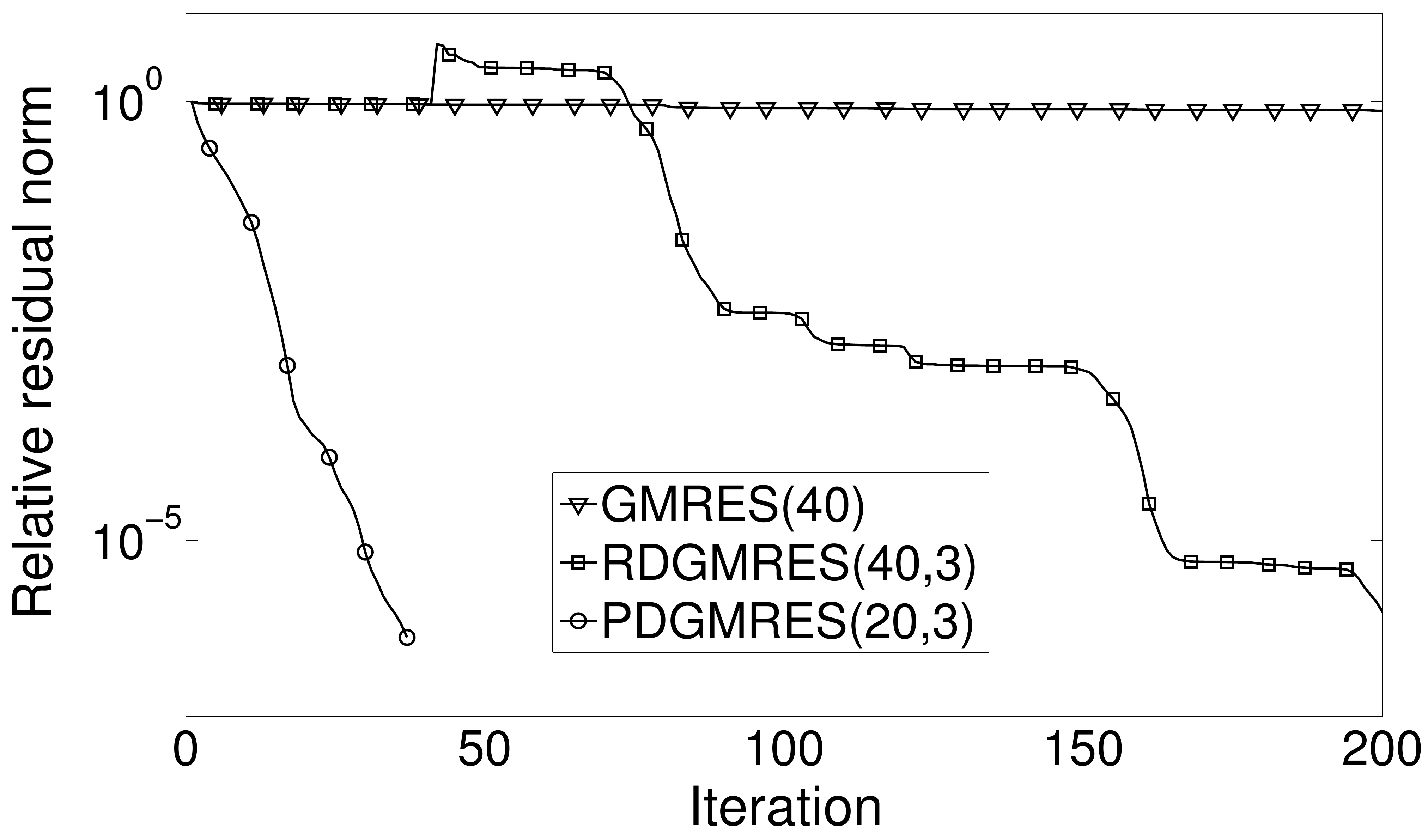}
\caption{Comparison of no deflation, harmonic Ritz deflation (RDGMRES) and (manual) physics-based deflation (PDGMRES) for SPE5 with modified permeability ($\sigma = 10^6$).}
\label{fig:Mandefl-SPE5}
\end{figure}

In figure, \ref{fig:Mandefl-SPE5}, we plot the convergence of the SPE5 case with modified permeability ($\sigma = 10^6$).
In line with the discussion of Figure \ref{fig:BO98-res530ritz}, $m$ is increased to $40$ in Figure \ref{fig:Mandefl-SPE5} to allow the harmonic Ritz vectors to converge to the eigenvectors corresponding to the extreme eigenvalues. RDGMRES reaches the tolerance after approximately 200 iterations. The convergence of PDGMRES is remarkably fast. This example illustrates that the span of the three manually constructed deflation vectors is a good approximation of the span of the eigenvectors corresponding to the extreme eigenvalues. Moreover, observe that we use $m = 20$ instead of $m = 40$. For $m = 40$, the convergence would be even better. We conclude that PDGMRES, compared to RDGMRES, not only achieves faster convergence but also does not impose any requirements on $m$, which, if $m$ can be lowered, decreases the computational cost of GMRES.

\subsubsection{Full simulation}
Next, we investigate the use of physics-based deflation in the full simulation. Model C (SAGD-SMALL) is used with the manually constructed deflation vectors from Figure \ref{fig:Mandefl-vecs}(c). Table \ref{tab:pdeflIX} summarizes the results for PDGMRES(30,3) and GMRES(30) with either a Jacobi or an AMG preconditioner for the pressure solve. We use a tolerance of $10^{-2}$ and the maximum number of iterations is set to two. The rationale behind the latter choice will be motivated in the next section, where we discuss the optimal PDGMRES settings.
\begin{table}[H]
	\centering
	\begin{tabular}{@{}llll@{}}
		\midrule[1.2pt]
		 & GMRES(30) & PDGMRES(30,3) & GMRES(30) \\
		\midrule
			Preconditioner		   & Jacobi & Jacobi & AMG \\
			Non-linears            & $279$ & $268$ & $291$\\
			Fails                  & $0$ & $0$ & $0$\\
			Outer linears          & $1,451$ & $1,085$ & $986$\\
			Inner linears          & $2,893$ & $2,164$ & $1,913$\\
			CPU time & $5.16$ & $6.58$ & $5.93$\\
		\midrule[1.2pt]
	\end{tabular}
	\caption{Comparison of GMRES(30) (with Jacobi preconditioner), PDGMRES(30,3) (with Jacobi preconditioner) and GMRES(30) (with AMG preconditioner) for the SAGD-SMALL case.}
	\label{tab:pdeflIX}
\end{table}
Each method uses two inner linear iterations in nearly all pressure solves. The size of the residual at the end of each pressure solve, however, will be smaller for PDGMRES and GMRES with AMG. This results in a reduction of the inner linear iterations for both methods. Note that the number of non-linear iterations varies. A better solution from the linear solve does reduce the non-linear iteration count in some cases, but our experiments are not conclusive.

In this case, both the time needed for the setup of AMG and the use of deflation offsets the iteration gain, resulting in higher CPU times. Note the difference with Table \ref{tab:deflIX}, where AMG was much faster. We see a better distinction (and better scalability for deflation) for larger cases. The case illustrates the potential of physics-based deflation as an alternative for AMG. Even though a significantly weaker preconditioner is used (Jacobi), PDGMRES still achieves a relatively good iteration reduction in the pressure solve. 

\subsubsection{Optimal PDGMRES settings for SAGD-SMALL case}
In this experiment, we investigate the optimal maximum number of inner linear iterations (denoted $N$). Commercial reservoir simulation typically uses only a single iteration to accelerate AMG in the pressure solve. To compete with this setup in terms of computational costs, the maximum number of iterations will have to be relatively low. The fact that physics-based deflation improves convergence from the start is especially advantageous (compared to harmonic Ritz deflation) for lower values of $N$. Table \ref{tab:maxcomp} compares GMRES without deflation to deflated GMRES with the 10 physics-based deflation vectors from Figure \ref{fig:Mandefl-vecs}(c) in the SAGD-SMALL case. A pressure solve tolerance of $10^{-2}$ is used.
\begin{table}[H]
	\centering
	\tiny{
	\begin{tabular}{@{}lllllllll@{}}
		\midrule[1.2pt]
		 & \multicolumn{2}{c}{Non-linears} & \multicolumn{2}{c}{Fails} & \multicolumn{2}{c}{Outer linears} & \multicolumn{2}{c}{Inner linears} \\
		\cmidrule(l){2-3}\cmidrule(l){4-5}\cmidrule(l){6-7}\cmidrule(l){8-9}
  		\addlinespace
		$N$ & PDGMRES & GMRES & PDGMRES & GMRES & PDGMRES & GMRES & PDGMRES & GMRES \\
		\midrule
		$1$ & $278$ & $272$ & $1$ & $1$ & $1,347$ & $1,568$ & $1,347$ & $1,568$ \\
		$2$ & $268$ & $279$ & $0$ & $1$ & $1,085$ & $1,451$ & $2,164$ & $2,893$ \\
		$3$ & $285$ & $279$ & $0$ & $0$ & $1,169$ & $1,383$ & $3,484$ & $4,119$ \\
		$5$ & $291$ & $287$ & $0$ & $0$ & $1,082$ & $1,297$ & $5,321$ & $6,394$ \\
		$10$ & $278$ & $277$ & $1$ & $1$ & $967$  & $1,221$ & $8,390$ & $11,623$ \\
		$20$ & $956$ & $974$ & $0$ & $0$ & $956$  & $974$ & $9,828$ & $15,721$ \\
		\midrule[1.2pt]
	\end{tabular}
	}
	\caption{PDGMRES and GMRES in the SAGD-SMALL case for varying $N$.}
	\label{tab:maxcomp}
\end{table}
As the maximum number of inner iterations is increased, the amount of outer linear iterations (roughly) decreases and the amount of inner linear iterations increases. If the pressure solve is allowed more iterations, then, in general, the residual will be smaller. As a result, the number of outer linears decreases. The number of non-linears stays approximately the same. Observe that the relative improvement of PDGMRES over GMRES  becomes more significant for higher values of $N$. For $N = 1$, the number of inner linears is decreased by $15\%$, compared to $37\%$ for $N = 20$. This can be explained by the fact that for $N = 10$, and especially for $N = 20$, PDGMRES often converges to the tolerance before reaching the maximum number of iterations, whereas GMRES does not. 

We aim to choose $N$ such that PDGMRES has a significant advantage over GMRES, while limiting the computational costs of the pressure solve. Although the advantage of deflation is more significant for $N = 10$ or $N = 20$, the additional pressure iterations increase the cost. In our experience, $2 \leq N \leq 5$ is sufficient to prevent failed non-linear iterations, while the number of inner linear iterations is kept relatively low. In general, we recommend taking $N = 2$. 

\subsubsection{SPE10 simulation}The complexity of modeling large models lies both in 
constructing $E^{-1}$ and solving of the deflated system $P_1 A\hat{x} = P_1 b$. 
Note that even when the deflation subspace dimension is relatively low, i.e. $d << n$, solving
the deflated system may still be difficult for challenging cases. On the other hand, increasing 
the deflation space dimension shifts the complexity to the construction of $E^{-1}$. However, it 
is possible to follow up the multiscale solution strategy in this case by ignoring the solution of 
the deflated system. This is the case where combing significantly weaker preconditioner such as Jacobi 
with the deflation method leads to a very poor performance. Hence, in this case, we combined deflation 
method with the AMG-preconditioned deflated system.

Although, this model has a clearly layered structure (see, figure \ref{fig:SPE10}) with two layers, the proposed 
physics-based deflation vectors for the SPE10 are equally distributed within the numerical domain using subdomain 
method (see, Figure \ref{fig:subdomain-levelset-deflation}(a)). The total number of deflation vectors considered 
in this study is (1) $d=8$, (2) $d=27$ and (3) $d=64$. The pressure system is solved approximately (as is common in 
practice) in any two-stage preconditioning methods. In this case, the pressure is also solved approximately.        

Table \ref{tab:pdeflIXSPE10} summarizes the results for PDGMRES(3,8), PDGMRES(3,27), PDGMRES(3,64) and CPR-AMG with an 
AMG preconditioner. As before, we use the termination criterion of either a tolerance of $10^{-2}$ or the maximum 
number of pressure solution iterations which is set to 3. Results and performance were compared with the current version 
of a modern reservoir simulator \cite{Schlumberger13} which uses CPR preconditioning with TRUE-IMPES decoupling, 1 V-cycle 
of AMG for solving the CPR pressure system and ILU(0) as a second-stage fully implicit (FIM) preconditioner \cite{Schlumberger13}. 
This reference solution strategy is referred to as CPR-AMG. The coarse grids in AMG solver are constructed by the parallel 
maximally independent set (PMIS) coarsening scheme with the Gauss-Seidel smoothing process. The coarse level is solved by 
FGMRES preconditioned by ILU(0) and maximum number of levels is limited to 50 by default settings.
The results for this case were obtained in serial runs using a desktop PC with 80GB RAM with Intel Xeon E5-2697 v2
2.70GHz CPU. Note that physics-based deflation is highly parallelizable, which is the topic of ongoing research.
\begin{table}[H]
	\centering
	\begin{tabular}{@{}lllll@{}}
		\midrule[1.2pt]
		                         & PDGMRES(3,8) & PDGMRES(3,27) & PDGMRES(3,64) & CPR-AMG \\
		\midrule
			Preconditioner		     & AMG           & AMG           & AMG            & AMG \\
			Non-linears            & $1,160$       & $1,159$       & $1,161$        & $1,159$\\
			Fails                  & $0$           & $0$           & $0$            & $0$\\
			Outer linears          & $12,986$      & $13,026$      & $13,029$       & $13,194$\\
			Inner linears          & $38,582$      & $38,775$      & $38,774$       & $39,248$\\
			Number of time steps   & $286$         & $286$         & $286$          & $286$\\ 
			CPU time               & $65,784.7$    & $66,035.3$    & $66,239.4$     & $56,916.1$\\
		\midrule[1.2pt]
	\end{tabular}
	\caption{Comparison of PDGMRES(3,8) (with AMG preconditioner), PDGMRES(3,27) (with AMG preconditioner), 
	                       PDGMRES(3,64) (with AMG preconditioner) and CPR-AMG 
	                       for the SPE10 case.}
	\label{tab:pdeflIXSPE10}
\end{table}

Each method uses approximately three inner linear iterations in all pressure solves. The size of the residual at the end of each pressure solve, however, will be smaller for PDGMRES with AMG. This results in a reduction of the inner linear iterations for all deflated methods. Note that the number of non-linear iterations in this case were the same.

In this case, both the time needed for the setup of AMG and the use of deflation offsets the iteration gain, resulting in higher CPU times. We expect 
the combination of multiscale solver and deflation method is expected to get a better solution strategy. Further analysis of optimal deflation 
vectors is required for this case and it is a subject of future research. In addition, it is important to note that the difference between the inner 
linear iterations for PDGMRES(3,27) and PDGMRES(3,64) is negligible leading to the conclusion that the dimension of the deflated space based on  
subdomain method must be relatively small for this case. This example also illustrates the potential of physics-based deflation vectors. The PDGMRES 
achieves in this case a relatively good iteration reduction in the pressure solve. 

\section{Conclusion}
We highlight the unfeasibility of the harmonic Ritz deflation method in practice, and demonstrate the potential of the physics-based levelset-subdomain algorithm. Deflation using harmonic Ritz vectors requires a full cycle of GMRES before the method can improve convergence. Furthermore, for more extreme eigenvalues computing accurate eigenvector approximations requires more iterations, which renders the harmonic Ritz deflation method less applicable precisely for those cases that would benefit the most from deflation. We show that the overhead of the method in our full simulations is larger than the time gained from the reduced number of iterations. In physics-based deflation, the deflation vectors are constructed before the start of the linear solve. Consequently, deflation can be used from the first iteration and onward. Moreover, the experiments show that the physics-based deflation vectors often perform better than the harmonic-Ritz deflation vectors. We conclude that physics-based deflation is, both in terms of feasibility and speed, better suited for commercial reservoir simulation than harmonic-Ritz deflation. Given the similarity with AMG when it comes to tackling harmful elements of the spectrum, as well as the good scalability properties, we believe that physics-based deflation combined with the robust parallel preconditioner for deflated system could potentially serve as an alternative for AMG. 

\section{Future Work}
In the present work, the physics-based deflation vectors are constructed manually. We show that the method yields excellent speedups compared to GMRES without deflation. The convergence improvements are the most evident for large, distinct, permeability jumps, such as in the SPE5 case with modified permeability and $\sigma > 10^4$. Experiments with the SAGD  and SPE10 cases show that physics-based deflation also works well for less well-defined permeability jumps, as most often encountered in real reservoirs. In commercial applications, manually assigning the deflation vectors is time consuming and will most likely lead to sub-optimal results. An automatic algorithm to apply the subdomain-levelset method, requiring minimal user input, will be introduced. We have highlighted the potential for a parallel implementation of this method in this paper, which will be further elaborated on. The combination of multiscale solver and deflation method is subject of future research.

\section*{Acknowledgments}
The authors would like to thank Schlumberger Abingdon Technology Center and Schlumberger-Doll Research 
for the support in establishing this work.

\bibliography{Deflation_paper}{}
\bibliographystyle{plain}

\end{document}